\documentclass[11pt]{article}
\usepackage{graphicx}
\usepackage{amsmath}
\usepackage{amsfonts}
\usepackage{epsfig}

\newcommand{\be}{\begin{equation}}
\newcommand{\ee}{\end{equation}}
\newcommand{\bes}{\begin{equation*}}
\newcommand{\ees}{\end{equation*}}
\newcommand{\beqn}{\begin{eqnarray}}
\newcommand{\eeqn}{\end{eqnarray}}
\newcommand{\beqns}{\begin{eqnarray*}}
\newcommand{\eeqns}{\end{eqnarray*}}

\newcommand{\fr}[1]{(\ref{#1})}

\newcommand{\lkr}{\left(}
\newcommand{\rkr}{\right)}
\newcommand{\lkv}{\left[}
\newcommand{\rkv}{\right]}
\newcommand{\lfi}{\left\{}
\newcommand{\rfi}{\right\}}

\newcommand{\EE}{{\bf E}}
\newcommand{\II}{{\bf 1}}

\newcommand{\Ni}{{\cal N}}

\newcommand{\sumk}{\sum_{k=0}^{2^j-1}}
\newcommand{\sumkp}{\sum_{k'=0}^{2^{j'}-1}}
\newcommand{\lamjeps}{\lambda_{j\varepsilon}}
\newcommand{\eps}{\varepsilon}
\newcommand{\tilf}{\tilde{f}}
\newcommand{\om}{\omega}
\newcommand{\tilom}{\tilde{\om}}
\newcommand{\chia}{\chi_{\eps, A}}

\newcommand{\card}{\mbox{card}}
\newcommand{\Var}{\mbox{Var}}

\newcommand {\bone}{\mathbf{e}}
\newcommand {\bu}{\mathbf{u}}
\newcommand {\bt}{\mathbf{t}}
\newcommand {\bx}{\mathbf{x}}
\newcommand {\bjp}{\mathbf{j'}}
\newcommand {\bkp}{\mathbf{k'}}
\newcommand {\bmp}{\mathbf{m'}}
\newcommand {\bJp}{\mathbf{J'}}
\newcommand {\bjpo}{\mathbf{j'_0}}
\newcommand {\bst}{\mathbf{s_2}}
\newcommand {\bspt}{\mathbf{s'_2}}
\newcommand {\bsstart}{\mathbf{s^*_2}}

\newcommand {\binfty}{\mbox{\mathversion{bold}$\infty$}}

\newcommand{\Kappa}{{\cal K}}

\newcommand{\bR}{{\mathbb R}}
\newcommand{\bN}{{\mathbb N}}
\newcommand{\bL}{{\mathbb L}}
\newtheorem{theorem}{Theorem}
\newtheorem{lemma}{Lemma}

\newtheorem{remark}{Remark}

\begin{document}

\title{\bf { Anisotropic Denoising in Functional Deconvolution Model with Dimension-free Convergence Rates}}

\author{{\em Rida Benhaddou and Marianna Pensky},  \\
         Department of Mathematics,
         University of Central Florida,  
          \\
         {\em Dominique Picard},\\
         University Paris-Diderot,
          CNRS, LPMA
         }

\date{}

\bibliographystyle{plain}
\maketitle

\begin{abstract}
 
In the present paper we consider the problem of    estimating   a periodic $(r+1)$-dimensional function  $f$ 
based on observations from its  noisy  convolution. We construct a wavelet 
estimator of $f$, derive minimax lower bounds for the $L^2$-risk when $f$   belongs to a  Besov 
ball of mixed smoothness and demonstrate that the wavelet estimator  is adaptive and asymptotically   near-optimal 
within a logarithmic factor, in a wide range of Besov balls. We prove in particular that choosing this 
type of mixed smoothness leads to rates of convergence which are free of the ''curse of dimensionality'' and, hence, are higher than
usual convergence rates when $r$ is large.

The problem studied in the paper is  motivated by seismic inversion which can be reduced to solution 
of noisy two-dimensional convolution equations that allow to draw inference on underground layer structures along the chosen profiles. 
The common practice in seismology  is to recover layer structures separately for each profile 
and then to combine the derived estimates into a two-dimensional function. By studying the two-dimensional version 
of the model, we demonstrate that this strategy usually  leads to estimators which are less accurate 
than  the ones obtained as two-dimensional functional deconvolutions.
Indeed, we show that   unless  the function $f$ is very smooth in the direction of the profiles, very
spatially inhomogeneous along the other direction and the number of profiles is very limited, 
the functional deconvolution solution has a much better precision compared  to a combination of $M$ solutions of separate  
convolution equations. A limited simulation study in the case of $r=1$ confirms theoretical claims of the paper.

\vspace{2mm} 


\vspace{2mm}{\bf AMS (2000) Subject Classification}: {Primary: 62G05,  Secondary: 62G08,62P35  }
\end{abstract}

\section{Introduction. }
\label{sec:introduction}
\setcounter{equation}{0}

Consider the problem of estimating a periodic $(r+1)$-dimensional function $f(\bu, x)$ with  
$\bu = (u_1, \cdots, u_r) \in [0,1]^r$  $x \in [0,1]$,  
based on observations from the following noisy convolution 
\be \label{eq1m}
y(\bu, t) =  \int^1_0 g(\bu, t-x) f(\bu, x)  dx + \varepsilon z(\bu, t),\ \ \ \bu \in [0,1]^r,  t \in [0,1].
\ee
Here,   $\eps$   is a positive small parameter such that asymptotically $\eps \to 0$, 
Function $g(. , .)$  in \fr{eq1m}  is assumed to be known and  $z(\bu,t)$ is  an  $r+1$-dimensional Gaussian white noise, i.e., a
generalized $r+1$-dimensional Gaussian field with covariance function
\begin{eqnarray*} 
\EE [z(\bu_1, t_1)z(\bu_2, t_1)] & = & \delta(t_1 - t_2)\ \prod_{l=1}^r \delta(u_{1l}- u_{2l}),
\end{eqnarray*}
where $\delta(\cdot)$ denotes the Dirac $\delta$-function and $\bu_{il} = (u_{i1}, \cdots, u_{ir}) \in [0,1]^r$, $i=1,2$.

Denote 
\bes
h(\bu, t) = \int^1_0  g(\bu, t-x)  f(\bu, x) dx.
\ees
Then, equation \fr{eq1m} can be rewritten as
\be \label{eq2}
y(\bu, t)= h(\bu, t) +  \varepsilon z(\bu, t)
\ee

In order to simplify the narrative, we start with the two dimensional version of equation \fr{eq1m}
\be \label{eq1}
y( u, t) = \int^1_0 g(u, t-x) f( u, x)  dx + \varepsilon z(u, t),\ \ \ u, t \in [0,1].
\ee
The sampling version of problem \fr{eq1} appears as 
\be \label{sampl}
y(u_l, t_i) = \int^1_0  g(u_l, t_i -x) f( u_l, x) dx  + \sigma z_{li}, \quad  l=1,\cdots, M,\ i=1,\cdots, N,
\ee
where $\sigma $ is a positive constant independent of $N$ and $M$,
$u_l = l/M$, $t_i = i/N$ and  $z_{li}$ are i.i.d normal variables with 
$\EE( z_{li})=0$, and $\EE( z_{l_1i_1}  z_{l_2i_2})=\delta(l_1-l_2)\delta(i_1-i_2)$.

Equation \fr{sampl} seems to be equivalent to $M$ separate  convolution equations
\be \label{separate}
y_l(t_i) = \int^1_0 f_l(x) g_l (t_i -x) dx  + \sigma z_{li}, \quad  l=1,\cdots, M,\ i=1,\cdots, N,
\ee 
with $y_l(t_i) = y(u_l, t_i)$, $f_l(x)= f(u_l, x)$ and $g_l (t_i -x) = g(u_l, t_i -x)$. 
This is, however, not true  since the solution of equation \fr{sampl} is a {\bf two-dimensional function} while 
solutions of equations \fr{separate} are $M$ unrelated functions $f_i(t)$.  In this sense, 
problem \fr{eq1} and its sampling equivalent \fr{sampl} are functional deconvolution problems.

Functional deconvolution problems have been introduced in Pensky and Sapatinas (2009) 
and further developed in  Pensky and Sapatinas (2010, 2011).  However, Pensky and Sapatinas (2009, 2010, 2011)
considered a different version of the problem where $f(u,t)$ was a function of one variable, i.e. 
$f(u,t) \equiv f(t)$. Their interpretation of functional deconvolution problem was motivated 
by solution of inverse problems in mathematical physics and multichannel deconvolution in engineering practices. 
Functional deconvolution problem of types \fr{eq1} and \fr{sampl} are motivated by experiments 
where one needs to recover a two-dimensional function using  observations  of its convolutions 
along profiles $u=u_i$. This situation occurs, for example, in geophysical explorations, in particular, 
the ones which rely on inversions of seismic signals (see, e.g., monographs of Robinson {\it et al.} (1996) and
Robinson (1999) and, e.g.,  papers of  Wason {\it et al.} (1984), Berkhout (1986)and  Heimer and Cohen (2008)).

In seismic exploration, a short duration seismic pulse is transmitted from the surface, reflected from
boundaries between underground  layers, and received by an array of sensors on the Earth surface. 
The signals are transmitted along straight lines called profiles. 
The received signals, called seismic traces, are analyzed to extract information about the underground
structure of the layers along the profile.  Subsequently, these traces can be modeled under simplifying assumptions as
noisy outcomes of convolutions between reflectivity sequences which describe configuration of the layers   
and the short  wave like function (called wavelet in geophysics) which corresponds to convolution kernel. The
objective of  seismic deconvolution is to estimate   the reflectivity sequences from the measured traces.
In the simple case of one layer and a single profile, the boundary will be described by an univariate function 
which is the solution of the convolution equation. The next step  is usually to combine the recovered functions 
which are defined on the set of parallel  planes passing through the profiles into a multivariate function which 
provides the exhaustive picture of the structure of the underground layers. 
This is usually accomplished by interpolation techniques. However, since the layers are intrinsically anisotropic 
(may have different structures in various directions) and spatially inhomogeneous (may experience, for example, sharp breaks),
the former  approach ignores the anisotropic and spatially inhomogeneous nature of the two-dimensional function describing the layer
and loses precision by analyzing each profile separately. 
\\ \\


The paper carries out  the following program:

\begin{itemize}
\item[i)] Construction of a feasible procedure $\widehat{f}(\bu, t)$ for estimating the $(r+1)$-dimensional function $f(\bu,t)$ 
which   achieves optimal rates of convergence (up to inessential logarithmic terms). 
 We require $\widehat{f}(\bu, t)$ to be adaptive with respect to smoothness constraints on $f$.
 In this sense, the paper is related to a multitude of papers which offered wavelet solutions to 
deconvolution problems (see, e.g.,  Donoho (1995), Abramovich and  Silverman (1998),
Pensky  and  Vidakovic (1999), Walter  and  Shen (1999),  Fan  and  Koo (2002),
Kalifa  and  Mallat (2003), Johnstone, Kerkyacharian, Picard  and 
Raimondo (2004), Donoho  and  Raimondo (2004), Johnstone   and  Raimondo
(2004), Neelamani, Choi  and  Baraniuk (2004) and Kerkyacharian, Picard
 and  Raimondo (2007)).

\item [ii)] Identification of the best achievable accuracy under smoothness constraints on $f$. 
We focus here on obtaining fast rates of convergence. In this context, we prove that considering 
multivariate functions with 'mixed' smoothness and hyperbolic wavelet bases allows to obtain rates 
which are free of dimension and, as a consequence, faster than the usual ones. In particular, the present 
paper is   related to anisotropic de-noising explored by, e.g., Kerkyacharian,   Lepski and   
Picard  (2001, 2008). We compare our functional classes as well as our rates with the results obtained there.

\item [iii)] Comparison of the two-dimensional version of the  functional deconvolution procedure studied in the present paper 
to the  separate solutions of convolution equations. We show especially that the former approach
delivers estimators with higher precision. For this purpose, in Section  \ref{sec:discrete_case},
 we consider a discrete version of functional deconvolution problem \fr{sampl} (rather than the continuous equation
\fr{eq1}) and compare its solution with solutions of  $M$  separate convolution equations \fr{separate}.
We show that, unless the function $f$ is very smooth in the direction of the profiles, very
spatially inhomogeneous along the other direction and the number of profiles is very limited, 
functional deconvolution solution has a better precision than the combination of $M$ solutions of separate  
convolution equations.

\end{itemize}

The rest of the paper is organized as follows.
In order to make the paper more readable and due to the application to seismic inversion, 
we start, in Section~\ref{sec:est_algorithm},  with the   two-dimensional version of the functional 
deconvolution problem \fr{eq1}, 
%
describe the construction of a  two-dimensional wavelet estimator of $f(u,t)$ given by equation \fr{eq1}.  
In Section~\ref{sec:lower_bounds}, we give a brief introduction on spaces of anisotropic smoothness. 
After that, we derive minimax lower bounds for the $L^2$-risk, based on observations from  
 \fr{eq1}, under the condition that  $f(u,t)$   belongs to a Besov ball of mixed regularity and
$g(u,x)$   has certain smoothness properties. In Section~\ref{sec:upper_bounds}, we prove that the hyperbolic wavelet estimator
derived in Section \ref{sec:est_algorithm} is adaptive and asymptotically
near-optimal within a logarithmic factor (in the minimax sense)  in a wide range of Besov balls.  
Section \ref{sec:discrete_case} is devoted to the discrete version of the problem \fr{sampl}
and comparison of functional deconvolution solution with the collection of individual deconvolution equations.
Section \ref{sec:multivariate} extends the results to the $(r+1)$-dimensional version of the problem \fr{eq1m}.
Section \ref{sec:simulations} contains a limited simulation study which supports theoretical claims of the paper.  
We conclude the paper by  discussion of the results in  Section~\ref{sec:discussion}.
Finally,  Section~\ref{sec:proofs} contains  the proofs of the theoretical results obtained
in the earlier sections.


\section{Estimation Algorithm. }
\label{sec:est_algorithm}
\setcounter{equation}{0}

 In what follows, $\langle \cdot,\cdot \rangle$ denotes the inner
product in the Hilbert space $L^2([0,1])$ (the space of
squared-integrable functions defined on the unit interval $[0,1]$),
i.e., $\langle f,g \rangle = \int_0^1 f(t)\overline{g(t)} dt$ for
$f,g \in L^2([0,1])$. We also denote the complex conjugate of $a$ by $\bar{a}$. 
Let $e_m(t) = e^{i 2 \pi mt}$ be a Fourier basis on  the interval  $[0,1]$. 
Let   $h_m(u) = \langle e_m, h(u, \cdot) \rangle$, 
$y_m(u) = \langle e_m, y(u, \cdot) \rangle$, $z_m(u)= \langle e_m, z(u, \cdot) \rangle$, 
$g_m(u)= \langle e_m, g(u, \cdot) \rangle$ and $f_m(u)= \langle e_m, f( u,\cdot) \rangle$ 
be functional Fourier coefficients of functions $h$, $y$, $z$, $g$ and $f$  respectively. 
Then, applying the Fourier transform   to  equation
\fr{eq2}, one obtains for any $u\in[0, 1]$  
\bes
y_m(u)= g_m(u)f_m(u) + \varepsilon z_m(u) 
\ees
and 
\be \label{hgf}
h_m(u)= g_m(u) f_m(u).
\ee

Consider   a bounded bandwidth periodized  wavelet basis (e.g., Meyer-type)  $\psi_{j, k}(t)$   
and finitely supported periodized  $s_0$-regular  wavelet  basis (e.g., Daubechies) $\eta_{j', k'}(u)$. 
The choice of the Meyer wavelet basis for $t$ is motivated by the fact that it allows easy evaluation of
the the wavelet coefficients in the Fourier domain while finitely supported wavelet basis gives more
flexibility in recovering a function which is spatially inhomogeneous in $u$.
Let $m_0$ and $m'_0$ be the lowest resolution levels for the two bases and denote the scaling functions 
for the bounded bandwidth wavelet  by $\psi_{m_0 -1,k} (t)$ and  the scaling functions 
for the finitely supported wavelet by $\eta_{m'_0 -1, k'} (u)$.
Then,    $f(u, x)$ can be expanded into wavelet series as  
\be \label{func1}
f(u, x) =  \sum_{j = m_0-1}^\infty \sum_{j' = m'_0 -1}^\infty   \sumk \sumkp\,  \beta_{j,k, j', k'} \psi_{j, k}(x) \eta_{j', k'}(u).
\ee
Denote $\beta_{j,k}(u) = \langle f(u,\cdot), \psi_{j, k} (\cdot) \rangle$, then, 
$\beta_{j, k, j', k'} = \langle \beta_{j,k}(\cdot), \eta_{j', k'}(\cdot) \rangle$.
If $\psi_{j, k, m} = \langle e_m,  \psi_{j, k} \rangle$ are Fourier coefficients of 
$\psi_{j, k}$, then, by formula \fr{hgf} and Plancherel's formula, one has
\be \label{betu}
\beta_{j, k}(u)=  \sum_{m \in W_j} f_m(u)\overline{ \psi_{j, k, m}} = \sum_{m \in W_j} \frac{h_m(u)}{g_m(u)}\ \overline{\psi_{j, k, m}},
\ee
where,   for any $j \geq j_0$, 
\be \label{eq:cj}
W_j = \lfi m: \psi_{jkm} \neq 0 \rfi \subseteq 2\pi/3 [-2^{j+2}, -2^j] \cup [2^j, 2^{j+2}],
\ee
 due to the fact that Meyer wavelets are band-limited (see, e.g., Johnstone, Kerkyacharian, Picard \& Raimondo
(2004), Section 3.1). 
Therefore, $\beta_{j,k, j', k'}$ are of the form
\be \label{beta}
\beta_{j,k, j', k'} = \sum_{m \in W_j} \overline{\psi_{j, k, m}} \ \int   \frac{h_m(u)}{g_m(u)}\  \eta_{j', k'}(u)du,
\ee
and allow the unbiased estimator 
\be \label{betaes}
\widetilde{\beta}_{j,k, j', k'} =  \sum_{m \in W_j} \overline{\psi_{j, k, m}} \ \int   \frac{y_m(u)}{g_m(u)}\  \eta_{j', k'}(u)du.
\ee
We now construct a hard thresholding estimator of $f(u,t)$ as 
\be   \label{Estres}
\widehat{f}(u, t) = \sum_{j = m_0-1}^{J-1} \sum_{j' = m'_0 -1}^{J'-1}  \sumk \sumkp\,  \widehat{\beta}_{jk, j'k'} \psi_{jk}(t) \eta_{j'k'}(u)
\ee
where 
\be \label{betathresh}
\widehat{\beta}_{j,k, j', k'}= \widetilde{\beta}_{j,k, j', k'} \II \left( \left| \widetilde{\beta}_{j,k, j', k'}\right| > \lambda_{j\varepsilon} \right). 
\ee
and the values of $J,J'$ and $\lamjeps$ will be defined later.

In what follows, we use the symbol $C$ for a generic positive
constant, independent of $\eps$, which may take different values at
different places.


\section{Smoothness classes and minimax lower bounds  }
\label{sec:lower_bounds} 

\subsection{Smoothness classes}
\label{subsec:functional_classes}
\setcounter{equation}{0}

It is natural to consider {\it anisotropic} multivariate functions, i.e., functions whose smoothness is  different in different directions. 
It is, however,  much more difficult to construct  appropriate spaces of mixed regularity which are meaningful for applications.
One of the objectives of the present paper is to prove that classes of mixed regularity  allow to obtain rates of convergence which are free of dimension.
This is specifically due to the application of hyperbolic wavelets, i.e., wavelets which allow  different resolution levels for each direction 
(see, e.g., Heping (2004)).

Although comprehensive study of functional classes of mixed regularity is not the purpose 
of this paper, below we provide a short introduction of functional classes that we are going to consider. 
Due to relation of this paper to anisotropic de-noising explored by   Kerkyacharian,   Lepski and   
Picard  (2001, 2008), we also  compare  classes of mixed regularity used therein to the 
Nikolski classes considered   in the papers  cited above.

 First, let us recall  definition of the Nikolski classes
$\Ni^{(s_1,\ldots,s_d)}_{(p_1,\ldots,p_d),\infty}$
(see  Nikolskii  (1975)).
In this section we consider $d$ dimensional multivariate functions. In what follows, we set $d=r+1$ or $d=2$.

Let $f$ be a measurable function defined on $\bR^d.$ For
any $x, y \in  \bR^d ,$ we define 
$$
\Delta_y f(x)= f(x + y) - f(x) .
$$
 If $l \in \bN $ then $\Delta_y^l$ is the $l-$iterated version of the
operator $\Delta_y .$ (Of course $\Delta_y^0 = I_d$  where $I_d$ is the identity operator.) 
Then,  Nikolski  classes can be defined as follows:\\
(recall that $\|g\|_{\bL^{p}(\bR^d,dx)}=\|g\|_{p}$ denotes $\left[\int_{\bR^d} |g(x_1,\ldots,x_d)|^p 
dx_1\ldots dx_d\right]^{1/p}$ for $1\le p <\infty$, with the usual modification for $p=\infty$.)

\begin{enumerate}   
\item 
Let $e_1,....e_d $ be the canonical basis of
$\bR^{d}$. For
$0 < s_i <\infty ; 1 \leq p_i \leq \infty $,
we say that $f$ belongs to
$\Ni^{s_i}_{p_i,\infty}$ if and only if
 there exists $ l \in \bN ,~ s_i < l, \hbox{ and } C(s_i,l) < \infty ,$ such that  
for any $h \in \bR $ one has
$$
\|\Delta_{he_i}^l f\|_{\bL^{p_i}(\bR^d,dx)}
\leq C(s_i,l)|h|^{s_i} .
$$

\item
$\Ni^{(s_1, \ldots,s_d)}_{(p_1,\ldots,p_d),\infty}=\cap_{i=1}^d  \Ni^{s_i}_{p_i,\infty}$ 
 \end{enumerate}

The  Nikolski classes defined above were investigated by  Kerkyacharian,   Lepski and   
Picard  (2001, 2008), they  are anisotropic but do not involve mixed smoothness. 
Quite differently,  in the present paper we shall  consider classes of mixed regularity defined as follows.
Denote $h=(h_1,\ldots,h_d)$, $t=(t_1,\ldots,t_d)$,  $s=(s_1,\ldots,s_d)$ and let $t_i>0$, $s_i >0$, $i=1, \cdots,d$.
For a subset $e\subset\{1,\ldots, d\}$, we set   $h^e$ to be the vector with coordinates $h_i$ when $i$ belongs to $e$, and 0 otherwise.
For   a fixed integer $l$ and  $1\le p\le \infty$, we denote
$$
\Delta_{h^e}^{l,e}f(x):=\left(\prod_{j\in e}\Delta_{h_je_j}^l\right)f(x), \quad 
\Omega^{l,e}(f,t^{e})_p:= \sup_{|h_j|\le t_j}\|\Delta_{h^e}^{l,e}f\|_p.
$$
Now, in order to construct Besov classes of mixed regularity, we choose $l \ge \max_j s_j$ and define 
\be \label{eq:Bes_infty} 
B^{s_1,\ldots,s_d}_{p,\infty} = \left\{ f\in \bL_p,\  \sum_{e\subset\{1,\ldots, d\}} \sup_{t>0}\, 
\sup_{j\in e}t_j^{-s_j}\Omega^{l,e}(f,t^{e})_p<\infty \right\}.
\ee
It is proved in, e.g., Heping (2004) that under appropriate (regularity) conditions which we are omitting here, 
 classes \fr{eq:Bes_infty} can be expressed in terms of hyperbolic-wavelet coefficients, thus, providing   
 a convenient  generalization of  the one-dimensional Besov $B^s_{p, \infty}$ spaces.
Furthermore,  Heping (2004) considers more general Besov classes 
of mixed regularity $B^{s_1,\ldots,s_d}_{p, q}$  that correspond to $q<\infty$ rather than $q=\infty$.
In this paper, we shall assume that the hyperbolic wavelet basis satisfies required regularity conditions
and follow  Heping (2004) definition of Besov spaces of mixed regularity 
\be \label{besov1}
B^{s_1,\ldots,s_d}_{p, q}= \left\{ f  \in L^2(U):\ \left( \sum_{j_1,\ldots,j_d} 2^{(\sum_{i=1}^d j_i[s_i+\frac12-\frac1p])
q} \left( \sum_{k_1,\ldots,k_d} 
\left| \beta_{j_1,k_1\ldots,j_dk_d}\right|^p \right)^{\frac{q}{p}} \right)^{1/q} <\infty\right\}.
\ee
Besov classes \fr{besov1} compare quite easily to the Nikolski classes: it is easy to prove that the former 
form a subset of the latter.

\subsection{Lower bounds for the risk:two-dimensional case}

Denote  $U=[0,1]\times [0,1]$ and
\be \label{sdef} 
s^*_i= s_i + 1/2 -1/p, \quad s'_i= s_i + 1/2 -1/p',\ \ \ i=1,2, \quad p'=\min\{ p, 2\}.
\ee
In what follows, we  assume that the function $f(u,t)$ belongs to a two-dimensional   
Besov ball as described above ($d=2$), so that   wavelet coefficients $ \beta_{j,k, j'k'}$  
satisfy the following condition
\be \label{assum1}
B^{s_1,s_2}_{p, q}(A)= \left\{ f  \in L^2(U):\ \left( \sum_{j, j'} 2^{(j s_1^* + j's_2^*)q} \left( \sum_{k, k'} 
\left| \beta_{j,k, j'k'}\right|^p \right)^{\frac{q}{p}} \right)^{1/q} \leq A \right\}.
\ee
Below, we construct  minimax lower bounds for the $L^2$-risk. For this purpose, we define the minimax $L^2$-risk
over the set $V$ as
$$
R_\eps (V) = \inf_{\tilde{f}} \, \sup_{f \in V} \EE \|
\tilde{f}  - f \|^2,
$$
where $\| g \|$ is the $L^2$-norm of a function $g(\cdot)$ and the
infimum is taken over all possible estimators $\tilde{f}(\cdot)$
(measurable functions taking their values in a set containing
$V$) of $f(\cdot)$.

Assume that functional Fourier coefficients $g_m(u)$ of function $g(u,t)$ are uniformly bounded from above and below, that is,  
there exist  positive constants $\nu$, and $C_1$ and $C_2$, independent of $m$ and $u$ such that  
\be \label{regsmo}
C_1 \left| m\right|^{-2\nu} \leq \left| g_m(u)\right|^2 \leq C_2 \left| m\right|^{-2\nu}.
\ee
Then, the following  theorem gives the minimax lower bounds for the $L^2$-risk of any estimator $\tilde{f}_n$ of $f$.

\begin{theorem}  \label{lowebd}
Let  $\min\{ s_1, s_2 \} \geq \max\{ 1/p, 1/2 \}$  with $ 1 \leq p,q \leq \infty$, let $A > 0$ and $s'_i$, $i=1,2$, be defined in \fr{sdef}. 
Then, under assumption \fr{regsmo}, as $\eps \rightarrow 0$
\be \label{eq:loweb}
R_\eps (B^{s_1,s_2}_{p, q}(A))\geq C A^2 \lkr \frac{\eps^2}{A^2} \rkr^d
\ee
where 
\be \label{d_value}
d = \min \lkr \frac{2s_2}{2s_2 +1}, \frac{2s_1}{2s_1 + 2\nu + 1}, \frac{2s'_1}{2s'_1 +2\nu} \rkr.
\ee 
\end{theorem}

Note that the value of $d$ in \fr{d_value} can be re-written as 
\be \label{d_val_cases}
d =  
\left\{ \begin{array}{ll}
\frac{2s_2}{2s_2 +1}, & \mbox{if}\ \ s_1 > s_2 (2\nu +1), \\
\frac{2s_1}{2s_1 + 2\nu + 1}, & \mbox{if}\ \ (\frac{1}{p}-\frac{1}{2})(2\nu +1) \leq  s_1 \leq s_2 (2\nu +1), \\
\frac{2s'_1}{2s'_1 +2\nu}, & \mbox{if}\ \ s_1 < (\frac{1}{p}-\frac{1}{2})(2\nu +1).
\end{array} \right.
\ee 

\begin{remark}
{\rm
Note  that the rates obtained here are in fact the worst rate associated to the one dimensional problem 
in each direction, which is not surprising since a function of only  one variable and constant in the other direction,  
e.g.,  $f(u_1,u_2)=h(u_1)$ belongs to $B^{s_1,s_2}_{p, q}(A)$ as soon as $h$ belongs to a ball of the usual 
one-dimensional Besov space $B^{s_1}_{p, q}$, for any $s_2$.
\\

Also it is worthwhile to observe that the third rate (involving $s'_1$) corresponds in dimension one to a ``sparse'' rate.  
Hence we observe here the so-called ``elbow phenomenon'' occurring only along the direction 2, because we are considering 
an $L^2$-loss and the problem has a degree of ill-posedness $\nu$ precisely in this direction.}
\end{remark}

\section{Minimax upper bounds. }
\label{sec:upper_bounds}
\setcounter{equation}{0}

Before deriving expressions for the minimax upper bounds for the risk, we formulate several useful 
lemmas which give some insight into the choice of the thresholds
$\lamjeps$ and upper limits $J$ and $J'$ in the sums in \fr{Estres}.

\begin{lemma} \label{lem:coef_var}
 Let $\widetilde{\beta}_{j,k, j', k'}$ be defined in \fr{betaes}. Then, under assumption \fr{regsmo}, one has 
\be \label{Varasym}
\Var \left(\widetilde{\beta}_{j,k, j', k'} \right)\asymp  \varepsilon^2  2^{2j \nu}.
\ee
\end{lemma}

\noindent
Lemma \ref{lem:coef_var} suggests that thresholds $\lamjeps$ should be chosen as 
\begin{eqnarray}  \label{tres}
\lambda_{j \varepsilon}= C_{\beta} \sqrt{ \ln ({1}/{\varepsilon})} \ 2^{j\nu}  \, \varepsilon
\end{eqnarray}
where $C_\beta$ is some positive constant independent of $\eps$. We choose $J$ and $J'$ as  
\be  \label{JJ}
2^{J}= (\varepsilon^2)^{-\frac{1}{2\nu +1}}, \quad 2^{J'}= (\varepsilon^2)^{-1}.  
\ee
Note that the choices of $J$, $J'$ and $\lambda_{j \varepsilon}$ are independent of the parameters, 
$s_1$, $s_2$, $p$, $q$ and $A$ of the Besov ball $B^{s_1s_2}_{p, q}(A)$, and therefore our estimator  
\fr{Estres} is adaptive with respect to those parameters.

The next two lemmas provide upper bounds for the wavelet coefficients and the large deviation inequalities 
for their estimators.

\begin{lemma} \label{lemA1}
Under assumption \fr{assum1}, one has
\bes
\sumk \sumkp\,   \left| \beta_{j,k, j', k'}   \right|^2 \leq A^{{2}} 2^{-2  (js'_1 + j' s'_2)}
\ees
for any $j, j'\geq 0$.
\end{lemma}

\begin{lemma} \label{lardiv}
Let $\widetilde{\beta}_{j,k, j', k'}$ and $\lambda_{j \varepsilon}$ be defined by formulae \fr{betaes} and  \fr{tres}, respectively. 
For some positive constant $\alpha$, define the set 
\be \label{Thetaset}
\Theta_{j,k, j'k', \alpha}= \{ \Theta:  \left| \widetilde{\beta}_{j,k, j', k'} - \beta_{j,k, j', k'}\right| > \alpha \lambda_{j\varepsilon} \}.
\ee
Then, under assumption \fr{regsmo}, as $\varepsilon \to 0$, one has
\be
\Pr \left( \Theta_{j,k, j'k', \alpha} \right)=  
O \left( \varepsilon^{\frac{\alpha^2 C^2_{\beta}}{2 \sigma_0^2}} \left[ \ln ({1}/{\varepsilon})\right]^{-\frac{1}{2}} \right)
\ee
where $ \sigma_0^2= \left( \frac{8 \pi}{3} \right)^{2 \nu} \frac{1}{C_1}$ and $C_1$ is defined in \fr{regsmo}.
\end{lemma}


Using the statements above, we can derive  upper bounds for the minimax  risk of the estimator \fr{Estres}.

\begin{theorem}  \label{uppbd} 
Let $\widehat{f}(. , .)$ be the wavelet estimator defined in \fr{Estres}, 
with 
$J$ and $J'$ given by \fr{JJ}. 
Let condition \fr{regsmo} hold and 
$\min\{ s_1, s_2 \} \geq \max\{ 1/p, 1/2 \}$, with 
$ 1 \leq p, q \leq \infty$. If $C_\beta$ in \fr{tres} is such that
\be \label{Cbeta_cond}
C_\beta^2 \geq 80 (C_1)^{-1} (2\pi/3)^{2 \nu} 
\ee
where $C_1$ is defined in \fr{regsmo}, then,  as $\eps \rightarrow 0$,
\be \label{eq:upper_bound}
\sup_{f \in  B^{s_1, s_2}_{p, q}(A))} \EE\| \widehat{f} - f\|^2 
\leq 
C A^2\ \lkr \frac{\eps^2 \ln(1/\eps)}{A^2} \rkr^d \, \ln \lkr \frac{1}{\eps} \rkr^{d_1}
\ee
where $d$ is defined in \fr{d_value} and 
\be \label{d1_value}
d_1 = \II(s_1 = s_2(2 \nu +1) ) + \II(s_1 = (2 \nu + 1)(1/p - 1/2)).
\ee
\end{theorem}

\begin{remark}
{\rm
Looking at the previous results, we conclude that the rates obtained by the wavelet estimator 
defined in \fr{Estres} are optimal, in the minimax sense, up to   logarithmic factors. 
These factors  are  standard and coming from the thresholding procedure.
}
\end{remark}


\section{Sampling version of the equation  and comparison with separate deconvolution recoveries }
\label{sec:discrete_case}
\setcounter{equation}{0}

Consider now the sampling version \fr{sampl}  of the problem \fr{eq1}.  In this case, the estimators of wavelet coefficients 
$\beta_{j,k, j', k'}$ can be constructed as 
\be \label{betaesdisver}
\widetilde{\beta}_{j,k, j', k'} =\frac{1}{M}  \sum_{m \in W_j} \overline{\psi_{j, k, m}} \ 
\sum^M_{l=1}  \frac{y_m(u_l)}{g_m(u_l)}\  \eta_{j', k'}(u_l).
\ee
In practice, $\widetilde{\beta}_{j,k, j', k'}$ are obtained simply by applying discrete wavelet transform to vectors $y_m(\cdot)/g_m(\cdot)$.

For any two sequences $a_n$ and $b_n$, one says that $a_n \asymp b_n$ as $n \to \infty$ if 
$0 < C_1 < a_n/b_n < C_2 < \infty$ for some constants $C_1$ and $C_2$ independent of $n$.  
Recall that the continuous versions \fr{betaes}  of   estimators \fr{betaesdisver} 
have $\Var \left(\widetilde{\beta}_{j,k, j', k'} \right)\asymp  \varepsilon^2  2^{2j \nu}$ (see formula \fr{Varasym}).
In order to show that equation \fr{sampl} is the sampling version of \fr{eq1} with $\eps^2 = \sigma^2 /(MN)$, one needs to show that,
in the discrete case, $\Var \left(\widetilde{\beta}_{j,k, j', k'} \right)\asymp  \sigma^2 (MN)^{-1}  2^{2j \nu}$. 
This indeed is accomplished by the following Lemma. 


\begin{lemma} \label{lem:coef_var_dicr}
 Let $\widetilde{\beta}_{j,k, j', k'}$ be defined in \fr{betaesdisver}. Then, under assumption \fr{regsmo}, 
as $MN \to \infty$, one has 
\be \label{Varadiscr}
\Var \left(\widetilde{\beta}_{j,k, j', k'} \right)\asymp  \sigma^2 (MN)^{-1}  2^{2j \nu}.
\ee
\end{lemma}

Using   tools developed in  Pensky and Sapatinas (2009)  and Lemma \ref{lem:coef_var_dicr}, it is easy to 
formulate the lower and the upper bounds for convergence rates of the estimator \fr{Estres} with 
$\widehat{\beta}_{j,k, j'k'}$ given by \fr{betathresh} and the values of $\lamjeps$ and  $J,J'$    defined in
\fr{tres} and \fr{JJ}, respectively. In particular, we obtain the following statement.


\begin{theorem}  \label{th:discrete}
Let  $\min\{ s_1, s_2 \} \geq \max\{ 1/p, 1/2 \}$  with $ 1 \leq p,q \leq \infty$, let $A > 0$ and $s^*_i$ be defined in \fr{sdef}. 
Then,  under assumption \fr{regsmo}, as $M N \rightarrow \infty$, for some absolute constant $C>0$ one has
\be \label{lower_discr}
R_{(MN)} (B^{s_1, s_2}_{p, q}(A))\geq C (\sigma^2 (MN)^{-1})^d.
\ee
Moreover, if  $\widehat{f}(. , .)$ is the wavelet estimator defined in \fr{Estres}, 
$\min\{ s_1, s_2 \} \geq \max\{ 1/p, 1/2 \}$,  and  $J$ and $J'$ given by \fr{JJ}, 
then, under assumption \fr{regsmo}, as $MN \rightarrow \infty$,
\be\label{upper_discr}
\sup_{f \in  B^{s_1, s_2}_{p, q}(A))} \EE\| \widehat{f} - f\|^2 
\leq C (\sigma^2 (MN)^{-1}  \ln(MN))^d \, (\ln(MN))^{d_1}.
\ee
where $d$ and $d_1$ are defined in \fr{d_value} and \fr{d1_value}, respectively.
\end{theorem}

Now, let us compare the rates in Theorem \ref{th:discrete} with the rates obtained by recovering each 
deconvolution  $f_l(t)= f(u_l, t)$, $u_l=l/M$, $l=1, \cdots, M$, separately, using equations \fr{separate}. 
In order to do this, we need to determine in  which space functions $f_l(x)$ are contained. The following lemma 
provides the necessary  conclusion.


\begin{lemma}\label{lem:onedBesov} 
Let  $f\in B^{s_1,s_2}_{p, q}(A)$ with  $s_1  \geq \max\{ 1/p, 1/2 \}$, $s_2> \max\{ 1/p, 1/2 \}$ and   $ 1 \leq p,q \leq \infty$.
Then, for any $l=1, ...., M$, we have
$$
f_l(t)= f(u_l, t)  \in B^{s_1}_{p, q}(\tilde{A}).
$$
\end{lemma}

Using Lemma \ref{lem:onedBesov} and standard arguments (see, e.g., 
Johnstone, Kerkyacharian, Picard  and  Raimondo (2004)), we obtain for each $f_l$
\begin{eqnarray*}
\sup_{f_l \in B^{s_1}_{p, q}(\tilde{A})} \EE\| \tilde{f}_l - f_l \|^2 \asymp  \left\{ \begin{array}{ll}
C  N^{-\frac{2s_1}{2s_1 + 2\nu + 1}}, & \mbox{if}\ \ s_1 \geq (\frac{1}{p}-\frac{1}{2})(2\nu +1), \\
C N^{-\frac{2s'_1}{2s'_1 +2\nu}}, 
& \mbox{if}\ \ s_1 < (\frac{1}{p}-\frac{1}{2})(2\nu +1).
\end{array} \right.
\end{eqnarray*}

Now, consider estimator $\tilde{f}$ of $f$ with $\tilde{f} (u_l, t_i) = f_l (t_i)$.  
If $f_u = \partial f/\partial u$ and $f_{uu} = \partial^2 f/\partial u^2$ exist and uniformly bounded for
$u \in [0,1]$, then rectangle method for numerical integration yields
$$
\EE\| \tilde{f}  - f  \|^2 =  M^{-1}\, \sum^M_{l=1} \EE\| \tilde{f}_l - f_l \|^2 + R_M,
$$
where 
$$
R_M \leq (12 M^2)^{-1}\ \lkv \EE\| \tilde{f}_u  - f_u  \|^2  + \sqrt{\EE\| \tilde{f}  - f  \|^2 \ 
  \EE\| \tilde{f}_{uu}  - f_{uu}  \|^2} \rkv. 
$$
If $M$ is large enough, then $R_M = o \lkr \EE\| \tilde{f}  - f  \|^2 \rkr$ as $M \to \infty$ and
we derive
\be \label{discr_separate}
\EE\| \tilde{f}  - f  \|^2    \asymp \left\{ \begin{array}{ll}
C  N^{-\frac{2s_1}{2s_1 + 2\nu + 1}}, & \mbox{if}\ \ s_1 \geq (\frac{1}{p}-\frac{1}{2})(2\nu +1), \\
C N^{-\frac{2s'_1}{2s'_1 +2\nu}}, 
& \mbox{if}\ \ s_1 < (\frac{1}{p}-\frac{1}{2})(2\nu +1).
\end{array} \right.
\ee
 
By straightforward calculations, one can check that the only case when convergence rates of separate 
deconvolution recoveries can possibly be better than that of the simultaneous estimator is when
$s_1 > s_2(2\nu +1)$. In this case, $s_1 > (\frac{1}{p}-\frac{1}{2})(2\nu +1)$, so that comparing the rates,
by straightforward calculations we derive that simultaneous recovery delivers better precision than separate ones 
unless
\be \label{comparison}
\lim_{\stackrel{M \to \infty}{N \to \infty}} M N^{-\frac{s_1 - s_2(2\nu+1)}{s_2(2s_1 + 2\nu +1)}} <1, \quad s_1 > s_2(2\nu +1).
\ee
It is easy to see that relation \fr{comparison} holds only if $s_1$ is large, $s_2$ is small 
and $M$ is relatively small in comparison with $N$.


\section{Extension to the $\mathbf{(r+1)}$-dimensional case }
\label{sec:multivariate}
\setcounter{equation}{0}

In this section, we extend the results obtained above   to the 
$(r+1)$-dimensional version of the model \fr{eq1m}. In this case, 
expanding both sides of equation \fr{eq1m} over Fourier basis, as before, we obtain  
 for any $\bu \in [0, 1]^r$  
$$
y_m(\bu)= g_m(\bu)f_m(\bu) + \varepsilon z_m(\bu). 
$$
Construction of the estimator follows the path of the two-dimensional case. With  $\psi_{j, k}(t)$ 
and $\eta_{j', k'}(u)$ defined earlier, we consider vectors $\bjp = (j'_1, \cdots, j'_r)$,   
$\bkp = (k'_1, \cdots, k'_r)$, $\bmp  = (m'_1, \cdots, m'_r)$ and $\bJp = (J'_1, \cdots, J'_r)$, 
and   subsets $\Upsilon(\bmp,\bJp)$ and $\Kappa(\bjp)$ of the set of 
$r$-dimensional vectors with nonnegative integer components:
\bes
\Upsilon (\bmp,\bJp) = \{ \bjp: m'_l \leq j'_l \leq J'_l, \ l=1, \cdots, r \},\quad
\Kappa(\bjp) =  \{ \bkp:  0 \leq k'_l \leq j'_l -1,  \ l=1, \cdots, r \}.
\ees
If $\binfty$ is the $r$-dimensional vector with all components being $\infty$, one can expand $f(\bu, t)$ into wavelet series as
\be   \label{func1m}
f(\bu, t) = \sum_{j = m_0-1}^{\infty} \sumk \, \sum_{\bjp \in \Upsilon (\bmp, \binfty)}\  \sum_{\bkp \in \Kappa(\bjp)}\   
\beta_{j, k, \bjp, \bkp} \psi_{jk}(t) \prod_{l=1}^r \eta_{j'_l, k'_l}(u_l),
\ee
where coefficients $\beta_{j, k, \bjp, \bkp}$ are of the form 
\be \label{betam}
\beta_{j, k, \bjp, \bkp} = \sum_{m \in W_j} \overline{\psi_{j, k, m}} \ \int_{[0,1]^d}\    \frac{h_m(\bu)}{g_m(\bu)}\, 
 \prod_{l=1}^r [\eta_{j'_l, k'_l}(u_l)]\, d\bu,
\ee
the set $W_j$ is defined by formula \fr{eq:cj} and $h_m(\bu) = \langle (f*g)(\cdot, \bu),  e_m (\cdot) \rangle$. 
Similarly to the two-dimensional case, we estimate $f(\bu, t)$ by 
\be   \label{Estresm}
\widehat{f}(\bu, t) = \sum_{j = m_0-1}^{J-1} \sumk \, \sum_{\bjp \in \Upsilon (\bmp, \bJp)}  
\sum_{\bkp \in \Kappa(\bjp)}\   \,  \widehat{\beta}_{j, k,\bjp, \bkp}\, \psi_{jk}(t)  \prod_{l=1}^r  \eta_{j'_l, k'_l}(u_l)
\ee
with 
\be \label{betathreshm}
\widehat{\beta}_{j,k, \bjp, \bkp}= \widetilde{\beta}_{j,k,  \bjp, \bkp} 
\II \left( \left| \widetilde{\beta}_{j,k,  \bjp, \bkp}\right| > \lambda_{j,\eps} \right). 
\ee
Here 
\be \label{betaesm}
\widetilde{\beta}_{j,k, \bjp, \bkp} =  \sum_{m \in W_j} \overline{\psi_{j, k, m}} \ \int   \frac{y_m(\bu)}{g_m(\bu)}\  
\prod_{l=1}^r [\eta_{j'_l, k'_l}(u_l)] d\bu
\ee
are the unbiased estimators of $\beta_{j, k, \bjp, \bkp}$, $J$ is defined in \fr{JJ}, 
$J'_l$ are such that $2^{J'_l} =  \eps^{-2}$, $l=1, \cdots, r$,   and $\lambda_{j,\eps}$ is given by formula \fr{tres}.

Assume, as before, that functional Fourier coefficients $g_m(\bu)$ of function $g(\bu,t)$ are uniformly bounded from above and below
\be \label{regsmom}
C_1 \left| m\right|^{-2\nu} \leq \left| g_m(\bu)\right|^2 \leq C_2 \left| m\right|^{-2\nu} 
\ee
and that function $f(\bu,t)$ belongs to an $(r+1)$-dimensional   Besov ball.  
As described in section \ref{subsec:functional_classes} to define these Besov balls, we introduce the vector 
$\bst = (s_{21}, \cdots, s_{2r})$ and denote by $\bspt$ and $\bsstart$ vectors with components  $s'_{2l} = s_{2l} + 1/2 -1/p'$ and 
$s^*_{2l} = s_{2l} + 1/2 -1/p$, $l=1, \cdots, r$, respectively, where $p'=\min\{ p, 2\}$. 
If $s_0 \geq \max_l s_{2l}$, then the $(r+1)$-dimensional   Besov ball of radius $A$ is characterized by its  
wavelet coefficients $ \beta_{j,k, \bjp, \bkp}$  as follows (see, e.g. Heping (2004))
\be \label{Besovm}
B^{s_1, \bst}_{p, q}(A)= \left\{ f  \in L^2([0,1]^{r+1}):\ \left( \sum_{j, \bjp} 2^{[j s_1^* + \bjp^T \bsstart  ]q} 
\left( \sum_{k, \bkp} \left| \beta_{j,k, \bjp, \bkp }\right|^p \right)^{\frac{q}{p}} \right)^{1/q} \leq A \right\}.
\ee 
It is easy to show that, with the above assumptions, similarly to the two-dimensional case, as $\varepsilon \to 0$, one has   
\begin{align}  
& \Var \left(\widetilde{\beta}_{j,k, \bjp, \bkp} \right)\asymp  \varepsilon^2  2^{2j \nu},  \ \ \ 
 \sumk \sumkp\,   \left| \beta_{j,k, \bjp, \bkp}   \right|^2 \leq A^{{2}} 2^{-2  (js'_1 +  \bjp^T \bsstart)}, \label{multvar}\\
& \Pr \left( \left| \widetilde{\beta}_{j,k, \bjp, \bkp} - \beta_{j,k, \bjp, \bkp}\right| > \alpha \lambda_{j\eps}  \right)=  
O \left( \eps^{\frac{\alpha^2 C^2_{\beta}}{2 \sigma_0^2}} \left[ \ln (1/\eps)\right]^{-\frac{1}{2}} \right). \label{multconcentr}
\end{align}
The upper and the lower bounds for the risk are expressed via 
\be \label{mins2}
s_{2,0} = \min_{l=1, \cdots, r} s_{2,l} = s_{2, l_0},
\ee 
where $l_0 = \arg \min s_{2,l}$. In particular, the following statements hold.

\begin{theorem}  \label{lowebdmult}
Let  $\min\{ s_1, s_{2, l_0} \} \geq \max\{ 1/p, 1/2 \}$  with $ 1 \leq p,q \leq \infty$. 
Then, under assumption \fr{regsmom}, as $\eps \rightarrow 0$,
\be \label{eq:lowebm}
R_\eps (B^{s_1,\bst}_{p, q}(A))\geq C A^2 \lkr \frac{\eps^2}{A^2} \rkr^D
\ee
where 
\be \label{d_valuem}
D = \min \lkr \frac{2s_{2,0}}{2s_{2,0} +1}, \frac{2s_1}{2s_1 + 2\nu + 1}, \frac{2s'_1}{2s'_1 +2\nu} \rkr.
\ee 
or, 
\be \label{D_valcasesm}
D =  
\left\{ \begin{array}{ll}
\frac{2s_{2,0}}{2s_{2,0} +1}, & \mbox{if}\ \ s_1 > s_{2,0} (2\nu +1), \\
\frac{2s_1}{2s_1 + 2\nu + 1}, & \mbox{if}\ \ (\frac{1}{p}-\frac{1}{2})(2\nu +1) \leq  s_1 \leq s_{2,0} (2\nu +1), \\
\frac{2s'_1}{2s'_1 +2\nu}, & \mbox{if}\ \ s_1 < (\frac{1}{p}-\frac{1}{2})(2\nu +1).
\end{array} \right.
\ee 
\end{theorem}

\begin{theorem}  \label{uppbdm} 
Let $\widehat{f}(. , .)$ be the wavelet estimator defined in \fr{Estresm}, 
with 
$J$  defined in \fr{JJ}, 
$J'_l$   such that $2^{J'_l} = (\eps^2)^{-1}$, $l=1, \cdots, r$,   and $\lambda_{j,\eps}$ given by formula \fr{tres}.
Let condition \fr{regsmo} hold and 
$\min\{s_1, s_{2, 0} \} \geq \max\{ 1/p, 1/2 \}$, with 
$ 1 \leq p, q \leq \infty$. If $C_\beta$ in \fr{tres} satisfies condition \fr{Cbeta_cond},
then,  as $\eps \rightarrow 0$,
\be \label{eq:upper_boundm}
\sup_{f \in  B^{s_1, \bst}_{p, q}(A))} \EE\| \widehat{f} - f\|^2 
\leq 
C A^2\ \lkr  A^{-2} \, \eps^2 \ln(1/\eps)  \rkr^D \, \ln \lkr 1/\eps \rkr^{D_1}
\ee
where $D$ is defined in \fr{d_valuem} and 
\be \label{D1_valuem}
D_1 = \II(s_1 = s_{2,0} (2 \nu +1) ) + \II(s_1 = (2 \nu + 1)(1/p - 1/2)) + \sum_{l\neq l_0}  \II(s_{2,l} = s_{2,0}).
\ee
\end{theorem}

\begin{remark} 
{\rm Observe that convergence rates in Theorems \ref{lowebdmult} and \ref{uppbdm}
depend on $s_1$, $p$, $\nu$ and $\min_l  s_{2l}$ but not on the dimension $r$.

It could be also natural to ask what would   the corresponding results be if  $s_1$ itself was  multidimensional, that is,
if one considers the case of convolution in more than one direction  where
$$
h(\bu, \bt) = \int_{[0;1]^d}  g(\bu, \bt-\bx)  f(\bu, x) d\bx, \quad \bt \in [0;1]^d;\ \bu \in [0;1]^r.
$$
Although  this is beyond the scope of this paper, let us just mention that, as soon as one establishes upper bounds for the variances 
of the wavelet coefficients like \fr{multvar} as well as concentration inequalities for the wavelet coefficients estimators 
like in \fr{multconcentr}, one expects to obtain convergence rates similar to Theorems \ref{lowebdmult} and \ref{uppbdm}
with $s_1$ replaced with $\min_k s_{1k}$. 
}
\end{remark}


\section{Simulations. }
\label{sec:simulations}
\setcounter{equation}{0}

In order to investigate finite-sample performance of our estimator, we carried out a limited simulation study. 
 We used WaveLab package for Matlab  and carried out simulations using degree 3 Meyer wavelet and degree 6 Daubechies  wavelets. 
We generated data  using equation \fr{sampl} with kernel $q(u,t) =   0.5\,  \exp(- |t|\,(1 + (u - 0.5)^2))$,
various functions $f(u,t)$ and various values of $M$, $N$ and $\sigma$.  In particular, we used $N=512$ ,
$M=128$ or $M =256$, $\sigma = 0.5$ or $\sigma = 1.0$ and  $f(u,t) =  f_1(u) f_2(t) $ where 
$f_1(u)$ and $f_2(t)$ are standard test functions routinely used in testing signal processing techniques 
(see, e.g., introduced by Donoho \& Johnstone (1994)). In particularly, we utilize functions
{\tt blip}, {\tt bumps},   and {\tt quadratic} with {\tt quadratic} just being a quadratic function $(y-0.5)^2$
scaled to have a unit norm.  Note that, though $f(u,t)$ is a product of two dimensional functions, the method 
does not ``know'' this and, therefore, cannot  take advantage of this information.

\begin{figure}[ht]
\[\includegraphics[height=4.0cm]{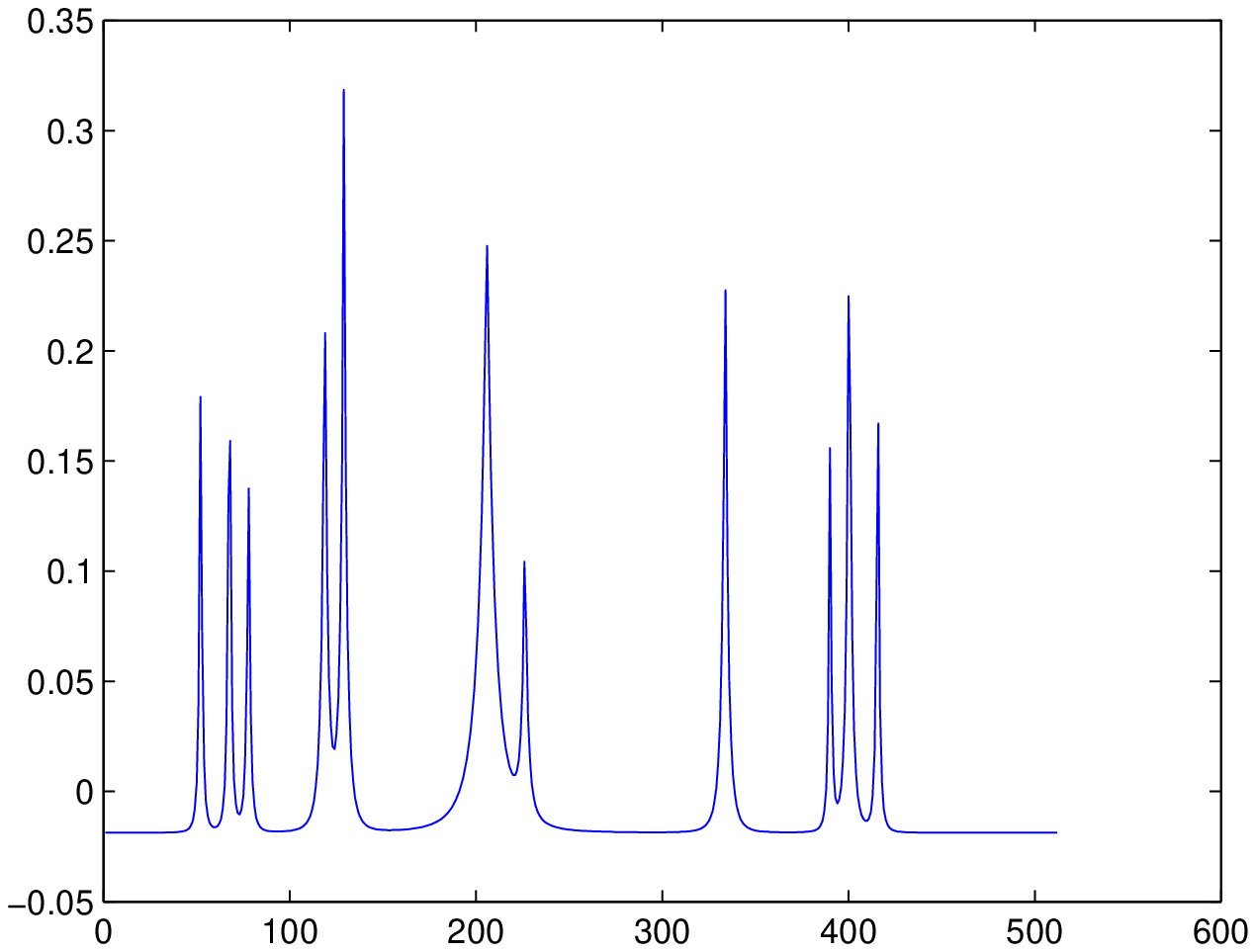} \hspace{3mm} \includegraphics[height=4.0cm]{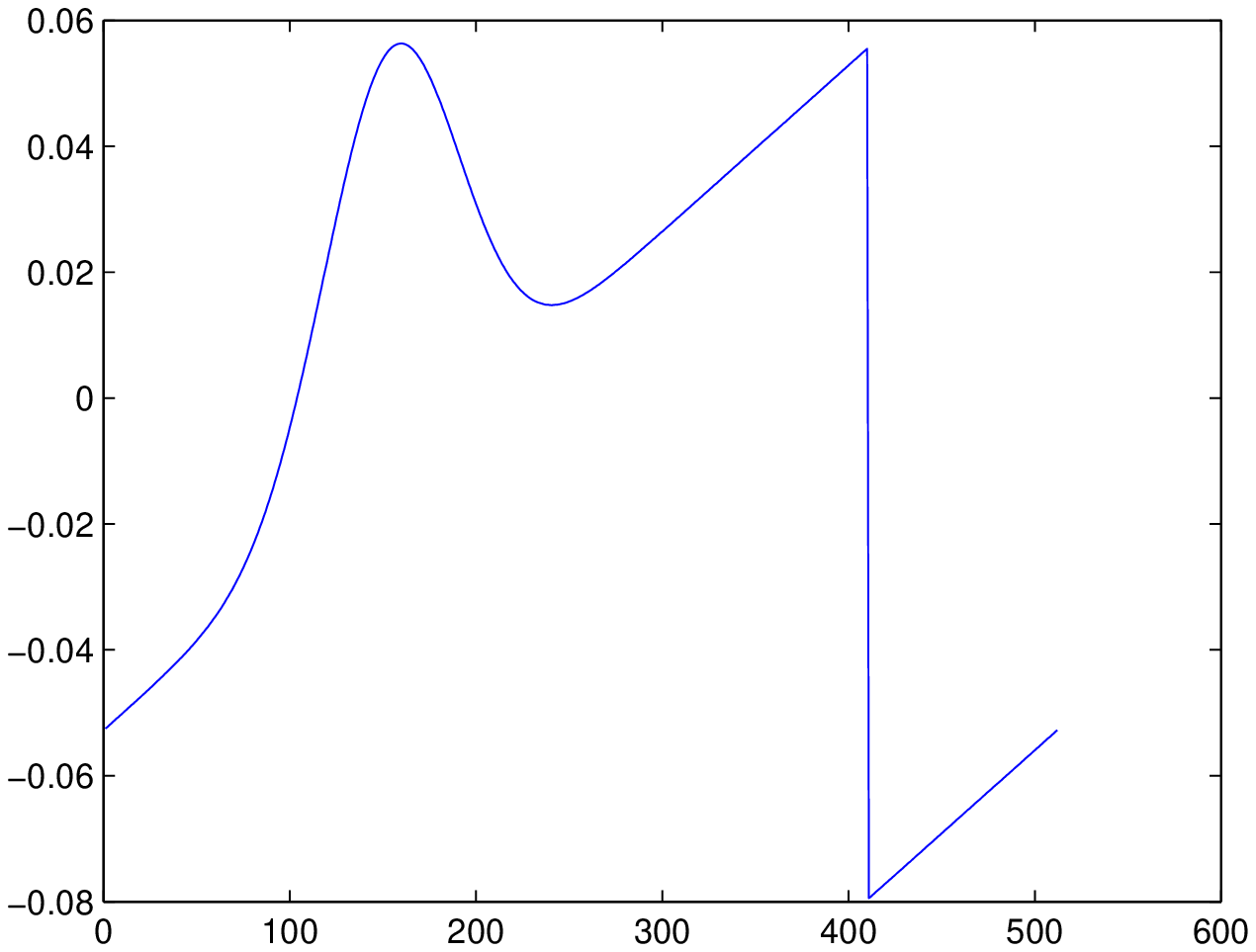}  
 \hspace{3mm} \includegraphics[height=4.0cm]{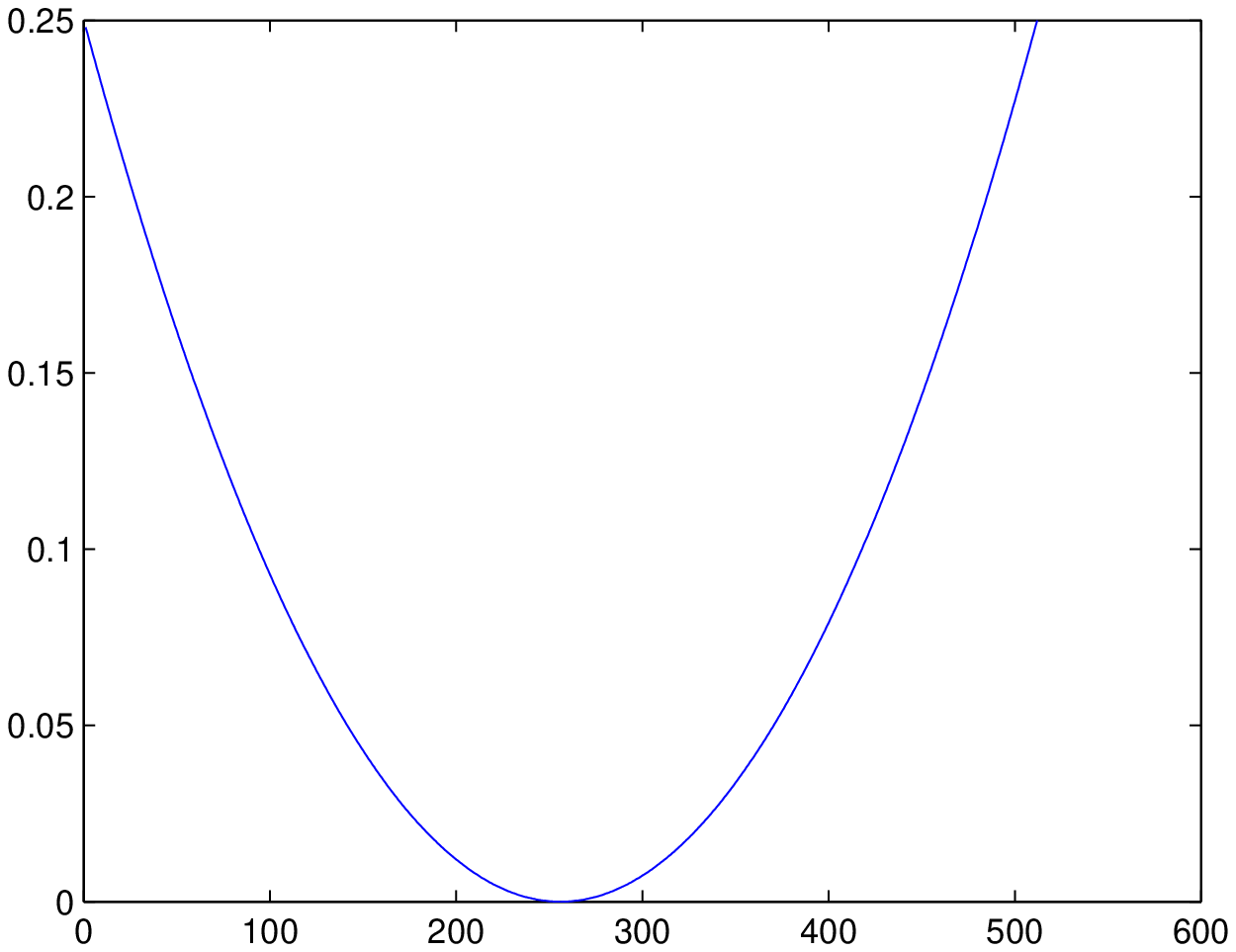} \]
%
\caption{ Test functions: {\tt bumps} (left), {\tt blip}  (middle),   
{\tt quadratic} (right) 
\label{fig1}}
\end{figure}

Graphs of all test functions are presented in Figure \ref{fig1}.

Table 1 contains simulations results. We generated data and constructed functional deconvolution estimator \fr{Estres}
and also $M$  Fourier-wavelet  deconvolution estimators of Johnstone, Kerkyacharian, Picard  and  Raimondo (2004)).
We evaluated mean integrated square error (MISE) $\EE \|\hat{f} - f\|^2$ of the functional deconvolution estimator and 
the average MISE of  $M$  Fourier-wavelet  deconvolution estimators.
Table 1 reports the averages of those errors over  100 simulation runs together with their standard deviations  
 (in the parentheses).  

Simulation results confirm that, as $M$ grows, functional deconvolution becomes more advantageous than $M$ separate 
deconvolutions. Indeed, while the error of a functional deconvolution estimator declines as $M$ grows, the average error of $M$ 
deconvolution estimators remains the same.


  \begin{table} [ht]
\begin{center}
\begin{tabular}{|c| c|c |c| c| c|}
  \multicolumn{6}{ c }{{\sc Functional deconvolution and $M$ separate deconvolutions   }}     \\
\hline
$M$ & $\sigma$ & MISE (functional) & MISE (separate) & MISE (functional) & MISE (separate) \\
\hline
 \multicolumn{2}{|c|}{ $N=512$} &
 \multicolumn{2}{|c|}{ $f_1 = $ {\tt Quadratic}, $f_2 = $ {\tt  Blip}} &
 \multicolumn{2}{|c|}{ $f_1 = $ {\tt Quadratic}, $f_2 = $ {\tt Bumps}}
 \\
\hline
 128 & 0.5  & 0.0535 (0.00148) &  0.0450 (0.00197) & 0.0534 (0.00123) &  0.0455 ( 0.00175) \\
\hline
128  & 1.0  &  0.213 (0.00614) &  0.181 (0.00816) &  0.212 (0.00589)  & 0.179 (0.00757) \\
\hline
256  & 0.5  & 0.0363 (0.00105) &  0.0452  (0.00148) & 0.0363 (0.000801) &  0.0451 (0.00133)   \\
\hline
256  & 1.0  &  0.145 (0.00331)  & 0.181 (0.00454)&  0.145 (0.00343)  &   0.180 (0.00458)  \\
\hline  \hline
\hline
 \multicolumn{2}{|c|}{$N=512$} &
 \multicolumn{2}{|c|}{ $f_1 = $ {\tt Blip}, $f_2 = $ {\tt  Blip}}&
 \multicolumn{2}{|c|}{ $f_1 = $ {\tt Blip}, $f_2 = $ {\tt Bumps}}
 \\
\hline
 128 & 0.5  &  0.0539 (0.00160) &   0.0453 (0.00190) & 0.0531 (0.00149) & 0.0447 (0.00208) \\
\hline
 128 & 1.0 &   0.214 (0.00695)  &   0.180 (0.00756)  & 0.214 (0.00661)  & 0.180 (0.00836) \\
\hline
256  & 0.5 &   0.0364 (0.000887) &  0.0452 (0.00120) & 0.0364 (0.00107) &  0.0452 (0.00149)\\
\hline
256  & 1.0 &   0.145 (0.00381)   &  0.180 (0.00572) &  0.145 (0.00420)  &  0.180 (0.00591)   \\
\hline \hline
 \multicolumn{2}{|c|}{  $N=512$} &
 \multicolumn{2}{|c|}{ $f_1 = $ {\tt Bumps}, $f_2 = $ {\tt  Blip}}&
 \multicolumn{2}{|c|}{ $f_1 = $ {\tt Bumps}, $f_2 = $ {\tt Bumps}}
 \\
\hline
 128 & 0.5  & 0.0535 (0.00145)  &  0.0452 (0.00144) & 0.0537 (0.00145) &  0.0454 (0.00197)  \\
\hline
 128 & 1.0 &  0.213  (0.00551)  &  0.179 (0.00727) & 0.214 (0.00683) &  0.181 (0.00751) \\
\hline
256  & 0.5 &  0.0363  (0.000925) & 0.0452 (0.00135) & 0.0364 (0.00101)  &  0.0451 (0.00144) \\
\hline
256  & 1.0 &   0.144 (0.00366)   &  0.180 (0.00467) & 0.146 (0.00355) & 0.181  (0.00479)  \\
\hline
\end{tabular}
\end{center}
\caption{ MISE   averaged  over $100$ runs. Third and fifth columns: average MISE 
of the functional deconvolution estimator. Fourth and sixth columns: average MISE of separate estimators for every $u$.
Standard deviations of the errors are listed in the parentheses. } \label{table1}
\end{table}


\section{Discussion. }
\label{sec:discussion}
\setcounter{equation}{0}

\begin{itemize}
\item[i)]
In the present paper, we constructed functional deconvolution estimators based 
on the hyperbolic wavelet thresholding procedure.  
We derived the lower and the upper bounds for the minimax convergence 
rates which confirm that estimators derived in the paper are   adaptive and asymptotically   near-optimal, 
within a logarithmic factor, in a wide range of Besov balls of mixed regularity.

\item[ii)] 
Although   results of Kerkyacharian,   Lepski and Picard  (2001, 2008)  have been obtained in a 
slightly different framework (no convolution), they can nevertheless be compared with the results presented above. 
Set $\nu=0$ to account for the absence of convolution, $p_i =p$  and  $d=r+1$.
Then, convergence rates in the latter can be identified as rates of a one-dimensional setting with a regularity parameter which is equal 
to the harmonic mean  
$$
\bar s = \left( \frac{1}{s_1}  + \cdots+ \frac{1}{s_d} \right)^{-1} < \min_{i=1, \cdots, d} s_i.
$$
In our case, the rates can also be identified as the rates in the one-dimensional setting with a regularity parameter  $\min_i s_i$
which is always larger than $\bar s$. Moreover, if $s_i =s$, one obtains $\bar s  = sd > s =  \min s_i$, showing that  
estimators of Kerkyacharian, Lepski and Picard  (2001, 2008)  in the Nikolski spaces are affected by ``the curse of dimensionality'' 
while the estimators  in the anisotropic Besov spaces of mixed regularity considered in this paper are free of 
dimension and, therefore, have higher convergence rates.

\item[iii)] The problem studied in the paper is  related to seismic inversion which can be reduced to solution 
of noisy convolution equations which deliver underground layer structures along the chosen profiles.  
The common practice in seismology, however, is to recover layer structures separately for each profile 
and then to combine them together. Usually, it is, however,  not the best strategy and leads 
to estimators which are inferior to the ones obtained as two-dimensional functional deconvolutions.
Indeed, as it is shown above,   unless function $f$ is very smooth in the direction of the profiles, very
spatially inhomogeneous along another dimension and the number of profiles is very limited, 
functional deconvolution solution has precision superior to combination of $M$ solutions of separate  
convolution equations. The precise condition when separate recoveries are preferable to the two-dimensional one 
is given by formula \fr{comparison} which,  essentially, is very reasonable. Really, if the number $M$ of profiles is small,
there is no reason to treat $f$ as a two-dimensional function.  Small value of  $s_2$  indicates that $f$ is very 
spatially inhomogeneous   and, therefore, the links between its values on different profiles are very weak. 
Finally, if $s_1$ is large, deconvolutions are quite precise, so that combination of various profiles 
cannot improve the precision.

\end{itemize}

\section*{Acknowledgments}

Marianna Pensky and Rida Benhaddou were partially supported by National Science Foundation
(NSF), grant  DMS-1106564.



\section{Proofs. }
\label{sec:proofs}
\setcounter{equation}{0}

\subsection{Proof of the lower bounds for the risk.}

In order to prove Theorem \ref{lowebd}, we consider two cases, the case
when $f(u,t)$ is dense in both variables (the dense-dense case) and the case 
when $f(u,t)$ is dense in $u$ and sparse in $t$ (the sparse-dense case).
The proof is based on Lemma A.1 of Bunea, Tsybakov and  Wegkamp~(2007) 
which we reformulate here for the case of squared risk.

\begin{lemma} \label{Tsybakov} [Bunea, Tsybakov, Wegkamp~(2007), Lemma A.1] 
Let   $\Omega$ be a set of functions of cardinality $\card(\Omega) \geq 2$ such that \\
(i) $\|f-g\|^2 \geq 4 \delta^2$,    \  \  \  for $f, g \in \Omega$, \ $f \neq g$,\\
(ii) the Kullback divergences $K( P_f, P_g)$ between the measures $P_f$ and $P_g$ satisfy the inequality
 $K( P_f, P_g) \leq   \log(\card(\Omega))/16$,    \  \  \ for \ \ $f, g \in \Omega$.\\
Then, for some absolute    positive constant $C$, one has 
 \begin{eqnarray*}
 \inf_{T_n} \sup_{f \in \ \Omega} \EE_f  \|T_n -f\|^2  \geq C \delta^2.
 \end{eqnarray*}
\end{lemma}

\underline{\bf  The dense-dense case.  } 
Let  $\omega$ be the matrix with components $\omega_{k, k'}=\{ 0,1\}$, $k=0,\cdots, 2^j-1$, $k' = 0,\cdots,2^{j'}-1$. 
Denote the set of all possible values $\omega$ by  $\Omega$ and let the functions $f_{j, j'}$ be of the form 
\begin{eqnarray} \label{denseden}
f_{jj'}(t, u)= \gamma_{jj'} \sum^{2^j-1}_{k=0} \sum^{2^{j'}-1}_{k'=0} \omega_{k, k'} \psi_{jk}(t) \eta_{j'k'}(u).
\end{eqnarray}
Note that matrix $\omega$ has $N =2^{j+j'}$ components, and, hence,  cardinality of the set of such matrices is $\card(\Omega)=2^N$. 
Since $f_{jj'} \in B^{s_1s_2}_{p, q}(A)$, direct  calculations show that $\gamma_{jj'} \leq A  2^{- j(s_1+1/2) - j'(s_2+1/2)}$, 
so that we choose   $\gamma_{jj'} = A  2^{- j(s_1+1/2) - j'(s_2+1/2)}$.
If $\tilde{f}_{jj'}$ is of the form \fr{denseden} with  $\tilde{\omega}_{k, k'} \in \Omega$ instead of $\omega_{k, k'}$,
then, the $L^2$-norm of the difference is of the form  
\begin{eqnarray*}
\| \tilde{f}_{jj'}-f_{jj'}\|^2  
&=& \gamma^2_{jj'} \sum^{2^j-1}_{k=0} \sum^{2^{j'}-1}_{k'=0} \II \left(  \tilde{ \omega}_{k, k'} \neq\omega_{k, k'}\right) 
= \gamma^2_{jj'} \rho( \tilde{ \omega}, \omega)
\end{eqnarray*}
where  $\rho( \tilde{ \omega}, \omega)=  \sum^{2^j-1}_{k=0} \sum^{2^{j'}-1}_{k'=0} \II \left(  \tilde{ \omega}_{k, k'} \neq\omega_{k, k'}\right)$ 
is the Hamming distance between the binary sequences $ \omega$ and $ \tilde{ \omega}$. In order to find a lower bound for the last expression, 
we apply the Varshamov-Gilbert lower  bound (see Tsybakov (2008), page 104) which states that one can choose 
 a subset $\Omega_1$ of $\Omega$, of cardinality at least $2^{N/8}$ such that $\rho( \tilde{ \omega}, \omega) \geq N/8$ 
for any $\omega , \tilde{ \omega} \in \Omega_1$. Hence, for any $\omega , \tilde{ \omega} \in \Omega_1$   one has
%
$\| \tilde{f}_{jj'}-f_{jj'}\|^2 \geq   \gamma^2_{jj'} 2^{j+j'}/8.$
Note that   Kullback divergence  can be written as 
\beqn  \label{Kulba1}
 K(f, \tilf) & = &  (2\eps^2)^{-1} \| (\tilf - f)*g \|^2.   
\eeqn
Since $|\om_{jj'} - \tilom_{jj'}| \leq 1$, plugging $f$   and $\tilf$ 
into \fr{Kulba1}, using Plancherel's formula and recalling that $|\psi_{j, k, m}| \leq 2^{-j/2}$, we derive 
\begin{eqnarray*}  
 K(f, \tilf)  & \leq & (2\eps^2)^{-1}  2^{-j}  \gamma^2_{jj'} \sum^{2^j-1}_{k=0} \sum^{2^{j'}-1}_{k'=0} 
\sum_{m \in W_j}\int_0^1 \eta^2_{j'k'}(u)\, g^2_m(u)\, du. 
 \end{eqnarray*}
Using \fr{regsmo}, we obtain
 \begin{eqnarray*}
 2^{-j} \sum_{m \in W_j}\int_0^1\eta^2_{j'k'}(u) g^2_m(u)du &\leq& C_2 2^{-j}\sum_{m \in W_j}|m|^{-2\nu} \int_0^1 \eta^2_{j'k'}(u)du 
 \leq  C_3  2^{-2\nu j},
 \end{eqnarray*}
so that
  \begin{eqnarray}   \label{kef}
 K(f, \tilf)  \leq   C \varepsilon^{-2}  \gamma^2_{jj'}  2^{j+j'} 2^{-2\nu j}.  
  \end{eqnarray}

Now, applying  Lemma \ref{Tsybakov} with 
\be \label{deltasq}
\delta^2=  \gamma^2_{jj'}{2^{j+j'}}/32 = A^2   2^{-  2s_1 j -  2s_2 j'}/32
\ee 
one obtains constraint $2^{-j(2s_1 +2\nu +1 )-j'(2s_2 +1)}  \leq   C \eps^2/A^2$ on
$j,j'$ and $\eps$ where $C$ is an absolute constant.  Denote 
\be \label{taue}
\tau_{\eps} = \log_2( C A^2 \varepsilon^{-2}).
\ee
Thus, we  need to choose combination of   $j$ and $j'$ 
which solves the following optimization problem 
\be  \label{opt1}
(j, j') = \arg\min \lfi    (2js_1+ 2j's_2)\quad \mbox{s.t.}\quad    j(2s_1 + 2\nu +1) + j'(2s_2+1) \geq \tau_{\eps}, \  j,j' \geq 0 \rfi.     
\ee
It is easy to check that solution of this linear constraint optimization problem is of the form
$\{ j, j'\}=\left\{ (2s_1 + 2 \nu +1)^{-1}  \tau_{\eps}, 0 \right\}$
if $s_2(2\nu +1) > s_1$, and $\{ j, j'\}=\left\{ 0, (2s_2 +1)^{-1} \tau_{\eps}  \right\}$ if $s_2(2\nu +1) \leq s_1$.
Plugging those values into \fr{deltasq}, obtain
\be \label{sol:dense}
\delta^2  = \left\{ \begin{array}{ll}
C A^2\,  (\eps^2/A^2)^{\frac{2s_2}{2s_2 +1}}, & \mbox{if}\ \ s_1 > s_2 (2\nu +1), \\
C  A^2\,  (\eps^2/A^2)^{\frac{2s_1}{2s_1 + 2\nu + 1}}, 
& \mbox{if}\ \    s_1 \leq s_2 (2\nu +1). 
\end{array} \right.
\ee

\underline{\bf  The sparse-dense case. } 
Let  $\omega$ be the vector with components $\omega_{ k'}=\{ 0,1\}$. Denote $\Omega$ the set of all possible $\omega$ and let the functions $f_{j, j'}$ be of the form
\begin{eqnarray}  \label{sparseden}
f_{jj'}(t, u)= \gamma_{jj'}  \sum^{2^{j'}-1}_{k'=0} \omega_{k'}\psi_{jk}(t) \eta_{j'k'}(u) 
\end{eqnarray}
Note that vector  $\omega$ has $N=2^{j'}$ components, and, hence, its cardinality is $\card(\Omega)=2^{N}$. 
Since $f_{jj'} \in B^{s_1s_2}_{p, q}(A)$, direct  calculations show that $\gamma_{jj'} \leq A  2^{- js^*_1 - j'(s_2+1/2)}$, so we choose $\gamma_{jj'} = A  2^{- js^*_1 - j'(s_2+1/2)}$.
If $\tilde{f}_{jj'}$ is of the form \fr{sparseden} with  $\tilde{\omega}_{k, k'} \in \Omega$ instead of $\omega_{k, k'}$,
then, calculating the $L^2$ norm of the difference similarly to dense-dense case, obtain   
\begin{eqnarray*}
\| \tilde{f}_{jj'}-f_{jj'}\|^2
&=& \gamma^2_{jj'} \sum^{2^{j'}-1}_{k'=0} \II \left(  \tilde{ \omega}_{ k'} \neq\omega_{ k'}\right) \geq 
 \gamma^2_{jj'}{2^{j'}}/8.
\end{eqnarray*}
Similarly  to dense-dense case,  using formulae \fr{regsmo} and \fr{Kulba1}, Plancherel's formula and    $|\psi_{j, k, m}| \leq 2^{-j/2}$,   derive
 \begin{eqnarray*}
 K(f, \tilf) \leq   (2\eps^2)^{-1}  \gamma^2_{jj'}   \sum^{2^{j'}-1}_{k'=0}2^{-j} \sum_{m \in W_j}\int_0^1\eta^2_{j'k'}(u) g^2_m(u)du
\leq  C  (2\eps^2)^{-1} \gamma^2_{jj'}  2^{j'} 2^{-2\nu j}.
  \end{eqnarray*}


Now, applying  Lemma \ref{Tsybakov} with 
\be \label{deltasq1}
\delta^2=  \gamma^2_{jj'}{2^{j'}}/32 = A^2   2^{-  2s'_1 j -  2s_2 j'}/32
\ee 
one obtains constraint $2^{-j(2s'_1 +2\nu)-j'(2s_2 +1)}  \leq   C \eps^2/A^2$ on
$j,j'$ and $\eps$ where $C$ is an absolute constant.  Thus, we  need to choose combination of   $j$ and $j'$ 
which delivers solution to the following linear optimization problem
%
\be     \label{opt2}
(j, j') = \arg\min \lfi (2js_1+ 2j's_2)\quad   \mbox{s.t.}\quad    j(2s'_1 + 2\nu) + j'(2s_2+1) \geq \tau_{\eps}, \  j,j' \geq 0\rfi.     
\ee
It is easy to check that solution of this linear constraint optimization problem is of the form
$\{ j, j'\}=\left\{ (2s'_1 + 2 \nu)^{-1}  \tau_{\eps}, 0 \right\}$
if $2 \nu s_2  > s'_1$, and $\{ j, j'\}=\left\{ 0, (2s_2 +1)^{-1} \tau_{\eps} \right\}$ if $2 \nu s_2 \leq s'_1$.
Plugging those values into \fr{deltasq1}, obtain
\be \label{sol:sparse}
\delta^2  = \left\{ \begin{array}{ll}
C A^2\,  (\eps^2/A^2)^{\frac{2s_2}{2s_2 +1}}, & \mbox{if}\ \ 2 \nu s_2 \leq s'_1, \\
C  A^2\,  (\eps^2/A^2)^{\frac{2s'_1}{2s'_1 + 2\nu}}, 
& \mbox{if}\ \    2 \nu s_2  > s'_1. 
\end{array} \right.
\ee
In order to complete the proof, recall expressions \fr{d_value} and \fr{d_val_cases} for $d$.


\subsection{Proofs of supplementary lemmas.}

{\bf Proof of Lemma \ref{lem:coef_var}. }
Let us derive an expression for the upper bound of the variance of \fr{betaes}. Subtracting  \fr{beta} from \fr{betaes} we obtain
\begin{eqnarray} \label{diffb}
\widetilde{\beta}_{j,k, j', k'}- \beta_{j,k, j', k'}= \varepsilon 
\sum_{m \in W_j} \overline{\psi_{j, k, m}} \ \int_0^1 \frac{z_m(u)}{g_m(u)}  \eta_{j', k'}(u)du. 
\end{eqnarray}
Now, before we proceed to the derivation of the upper bound of the variance, let us first state a result that will be used in our calculation. 
Recall from stochastic calculus that for any function $F(t, u) \in L^2 ([0,1]\times[0,1]) $, one has 
\begin{eqnarray} \label{exsto}
\EE \left[ \int_0^1 \int_0^1 F(t, u) dz(t, u)du\right]^2 &=& 
\int_0^1 \int_0^1 F^2(t, u)dt du. 
\end{eqnarray}
Hence, recalling that $z_m(u) = \int z(u,t) e_m(t) dt$, choosing  
$$
F(t, u)= \sum_{m \in W_j} \overline{\psi_{j, k, m}} \ \frac{e_m(t)}{g_m(u)} \eta_{j', k'}(u), 
$$ 
squaring both sides of \fr{diffb}, taking expectation and using the relation \fr{exsto}, we obtain
\begin{eqnarray*} 
\Var \left(\widetilde{\beta}_{j,k, j', k'} \right) 
&=& 
\varepsilon^2\ \EE \left| \sum_{m \in W_j} \overline{\psi_{j, k, m}} \int_0^1 \int_0^1 \frac{\eta_{j', k'}(u)}{g_m(u)} e_m(t) dz(u, t)du \right|^2\\
 &=& 
\varepsilon^2 \int_0^1 \int_0^1 \sum_m \sum_{m'} \frac{  \overline{\psi_{j, k, m}} \psi_{j, k, m'}}{g_m(u) \overline{g_{m'}} (u)} 
\overline{e_m}(t) e_{m'}(t) |\eta_{j', k'}(u)|^2 dt du\\
&=& 
\varepsilon^2 \sum_{m \in W_j} \left|  \psi_{j, k, m}\right|^2 \int_0^1 \frac{ |\eta_{j', k'}(u)|^2}{ |g_m(u)|^2}du,
\end{eqnarray*}
since in the double summation above, all terms involving  $m \neq m'$   vanish due to  $\int_0^1 e_m(t) e_{m'}(t) dt=0$. 
Consequently,
Taking into account \fr{eq:cj}, \fr{regsmo} and the fact that $|\psi_{j, k, m}| \leq 2^{-j/2}$, obtain
\be \label{varexpr}
\Var \left(\widetilde{\beta}_{j,k, j', k'} \right)  
\asymp  \varepsilon^2 \sum_{m \in W_j} |\psi_{j, k, m}|^2|m|^{2\nu} \int^{1}_{0} \left| \eta^2_{j', k'}(u) \right| du
   \asymp  \varepsilon^2  2^{2j \nu}
\ee
so that \fr{Varasym} holds.
\\


{\bf Proof of Lemma \ref{lemA1}  } First note that,  under assumption \fr{assum1}, one has  
  \begin{eqnarray*}
\sum_{k, k'} \left| \beta_{j,k, j', k'}   \right|^p \leq A^{p} 2^{-p \left[ (js_1 + j' s_2) + (\frac{1}{2}-\frac{1}{p})(j+ j') \right] }
 \end{eqnarray*}  
If $p\leq 2$, one has $p'=p$, $s'_i = s_i + 1/2-1/p$, $i=1,2$, and 
\begin{eqnarray*}
\sum_{k, k'} \left| \beta_{j,k, j', k'}   \right|^2 
& \leq & \sum_{k, k'} \left| \beta_{j,k, j', k'}   \right|^p \left \{ \max_{k, k'} \left| {\beta}_{j,k, j', k'}\right|^p \right\}^{(2-p)/p}
\leq A^{{2}} 2^{-2 (js'_1 + j' s'_2)}.
\end{eqnarray*}  
If  $p \geq 2$, then $p'=2$, $s'_i = s_i$, $i=1,2$, and,  applying the Cauchy-Schwarz inequality, one obtains
 \begin{eqnarray*}
 \sum_{k, k'} \left| \beta_{j,k, j', k'}   \right|^2 
&\leq & 
\left( \sum_{k, k'} \left| \beta_{j,k, j', k'}   \right|^p \right)^{2/p}\left( \sum_{k, k'} 1\right)^{(1-{2/p})}
\leq A^{2} 2^{-2 \left[ (js_1 + j' s_2) \right]},
\end{eqnarray*}  
which completes the proof.
\\


{\bf Proof of Lemma \ref{lardiv} } Observe that $\widetilde{\beta}_{j,k, j', k'}- \beta_{j,k, j', k'}$ 
is a zero-mean Gaussian random variable with variance given by \fr{varexpr}, so that 
\be  \label{varup}
\Var \left(\widetilde{\beta}_{j,k, j', k'} \right)  \leq 
\varepsilon^2\left( \frac{8 \pi}{3} \right)^{2 \nu} \frac{2^{2\nu j}}{C_1} = 
  \sigma_0^2  \varepsilon^2 2^{2\nu j} 
\ee
Denoting by $\bar{\Phi} (x) = 1 - \Phi(x)$ where $\Phi(x)$ is the standard normal c.d.f. 
and recalling that $\bar{\Phi} (x) \leq (x \sqrt{2 \pi})^{-1} \exp(-x^2/2)$ if $x>0$,
we derive
\begin{eqnarray*} 
\Pr \left( \Omega_{jk, j'k', \alpha} \right)
& = & 
\Pr \left( \left| \xi_{j, k, j', k'}\right| > \alpha \lambda_{j\varepsilon} \right) 
= 2  \bar{\Phi}\left( \alpha \lambda_{j\varepsilon} (\sigma_0   \varepsilon  2^{\nu j})^{-1} \right)\\
& \leq &
2  \bar{\Phi}\left(  \alpha C_\beta (\sigma_0)^{-1} \sqrt{ \ln (1/\eps)} \right)
\leq  \frac{2 \sigma_0}{\alpha C_{\beta}  \sqrt{2\pi \ln (1/\varepsilon)}}\ 
\eps^{ \frac{\alpha^2  C^2_{\beta}}{2\sigma^2_0} }    
\end{eqnarray*}  
which completes the proof. \\

\subsection{Proof of upper bounds for the risk.}

{\bf Proof of Theorem \ref{uppbd} } Denote
\be \label{chi_epsA}
\chia = A^{-2} \eps^2 \, \ln(1/\eps),
\ee
\be \label{j0j10}
2^{j_0} = (\chia)^{-\frac{d}{2 s'_1}},\quad 2^{j'_0} = (\chia)^{-\frac{d}{2 s'_2}}
\ee
and observe  that with $J$ and $J'$ given by \fr{JJ},   the estimation error can be decomposed into the sum of four components  as follows
\begin{eqnarray}  \label{eserto}
\EE \| \widehat{f}_n - f \|^2 &\leq& \sum_{j,k, j',k'} \EE \| \widehat{\beta}_{j,k, j', k'}- \beta_{j,k, j', k'}   \|^2  
  \leq   R_1 + R_2 + R_3 + R_4, 
\end{eqnarray}
where 
\begin{eqnarray*}
R_1&= &  \sum^{2^{m_0}-1}_{k=0} \sum^{2^{m'_0}-1}_{k'=0} \Var( \widetilde{\beta}_{m_0, k, m'_0, k'}), \\ 
R_2&=& \sum^{J-1}_{j= m_0 } \sum^{J'-1}_{j'=m'_0} \sum_{k, k'} \EE \left[ \left| \widetilde{\beta}_{j,k, j', k'}- \beta_{j,k, j', k'}    \right|^2 
\II \left( \left| \widetilde{\beta}_{j,k, j', k'}\right| > \lambda_{j\varepsilon}  \right)\right], \\
R_3 &=& \sum^{J-1}_{j=m_0 } \sum^{J'-1}_{j'=m'_0} \sum_{k, k'} \left| \beta_{j,k, j', k'}   \right|^2 
\Pr \left( \left| \widetilde{\beta}_{j,k, j', k'}\right| < \lambda_{j\varepsilon}    \right), \\
R_4 &=& \left( \sum^{\infty}_{j=J}\sum^{J'-1}_{j'=m'_0}  + \sum^{J-1}_{j=m_0}\sum^{\infty}_{j'= J'}+ 
\sum^{\infty}_{j= J}\sum^{\infty}_{j'= J'} \right)\sum_{k, k'} \left| \beta_{j,k, j', k'}   \right|^2.  
\end{eqnarray*}  
For  $R_1$, using \fr{Varasym}, derive, as $\eps \to 0$, 
\be \label{er1}
R_1  \leq C \varepsilon^2 = O \lkr A^2\, \chia^d \rkr.
\ee
To calculate $R_{4}$, we apply Lemma  \ref{lemA1}   and use \fr{JJ} obtaining, as $\eps \to 0$, 
\begin{eqnarray}
R_4&=& O \left( \left( \sum_{j\geq J} \sum_{j' \geq m'_0} + \sum_{j\geq m_0} \sum_{j' \geq J'} \right)A^2 2^{-2js'_1 - 2j's'_2} \right)
= O \left( A^22^{-2Js_1} +A^2 2^{-2J's_2} \right)\nonumber\\
&=& O \left( A^2 (\varepsilon^2)^{\frac{2s'_1}{2\nu +1}} +A^2 (\varepsilon^2)^{2s'_2}  \right) = O \lkr A^2 \chia^d \rkr.
\label{er4}
\end{eqnarray}
Then, our objective is to prove that, as $\eps \to 0$, one has $R_i = O \lkr A^2 \chia^d [\ln (1/\eps)]^{d_1} \rkr.$

Now, note that  each $R_2$ and $R_3$ can be partitioned into the sum of two errors as follows
\begin{eqnarray}  \label{Triaine}
R_2  \leq  R_{21} + R_{22}, \quad R_3 \leq   R_{31} + R_{32},
\end{eqnarray}
where
 \begin{eqnarray}
 R_{21}&=& \sum^{J-1}_{j=m_0} \sum^{J'-1}_{j'=m'_0} \sum_{k, k'} \EE \left[ \left| \widetilde{\beta}_{j,k, j', k'}- \beta_{j,k, j', k'}    \right|^2 
\II \left( \left| \widetilde{\beta}_{j,k, j', k'} - \beta_{j,k, j', k'}    \right| > \frac{\lambda_{j\varepsilon}}{2}  \right)\right] \ \ \  
\label{r21}\\
 R_{22}&=& \sum^{J-1}_{j=m_0} \sum^{J'-1}_{j'=m'_0} \sum_{k, k'} \EE \left[ \left| \widetilde{\beta}_{j,k, j', k'}- \beta_{j,k, j', k'}    \right|^2 
\II \left( \left|{\beta}_{j,k, j', k'}\right| >\frac{1}{2} \lambda_{j\varepsilon}  \right)\right].\ \ \  
\label{r22}\\
 R_{31}& = & \sum^{J-1}_{j=m_0} \sum^{J'-1}_{j'=m'_0} \sum_{k, k'} \left| \beta_{j,k, j', k'}   \right|^2 
\Pr \left( \left| \widetilde{\beta}_{j,k, j', k'} - \beta_{j,k, j', k'}\right| > \frac{\lambda_{j\varepsilon}}{2}    \right),\ \ \  
\label{r31}\\
 R_{32}&=& \sum^{J-1}_{j=m_0} \sum^{J'-1}_{j'=m'_0} \sum_{k, k'} \left| \beta_{j,k, j', k'}   \right|^2 
\II \left( \left| {\beta}_{j,k, j', k'}\right| \leq \frac{3 \lambda_{j\varepsilon}}{2}    \right).\ \ \  
\label{r32}
\end{eqnarray}  
 Combining \fr{r21} and \fr{r31}, and applying Cauchy-Schwarz inequality and Lemma \ref{lardiv} with $\alpha=1/2$, one derives
\begin{eqnarray*} 
R_{21} + R_{31} 
&=& 
O \left( \sum^{J-1}_{j=m_0}\sum^{J'-1}_{j'=m'_0}  2^{j+j'} \, \varepsilon^{\frac{C^2_{\beta}}{16 \sigma^2_0}} \left[ \ln ({1}/{\varepsilon})\right]^{-\frac{1}{4}} 
\sqrt{ \varepsilon^4 2^{4j\nu+j'}}     \right)\\
&=& 
O \left( 2^{J(2\nu+1)}\, 2^{3J'/2}\, (\varepsilon)^{2 + \frac{C^2_{\beta}}{16 \sigma^2_0}} \right) 
=  O \left(  (\varepsilon^{2})^{ \frac{C^2_{\beta}}{32 \sigma^2_0}-\frac{3}{2}}\right).  
\end{eqnarray*}
Hence, due to condition \fr{Cbeta_cond}, one has, as $\eps \to 0$, 
\be \label{er2}
R_{21} + R_{31} \leq C \varepsilon^2 = O \lkr A^2 \chia^d \rkr.
\ee
For the sum of  $R_{22}$ and  $R_{32}$, using \fr{Varasym} and \fr{tres}, we obtain
\be \label{delta_expr}
\Delta = R_{22}+ R_{32} = O \lkr  \sum^{J-1}_{j=m_0} \sum^{J'-1}_{j'=m'_0} \sum_{k, k'}  \min \left \{ \beta^2_{j,k, j', k'}, 
\eps^2 \ln ( {1}/{\varepsilon}) \ 2^{2j\nu} \right  \}   \rkr.
\ee
Then, $\Delta$ can be partitioned  into the sum of  three components  $\Delta_{1}$,  $\Delta_{2}$ and $\Delta_3$  
according to three different sets of indices:
\begin{eqnarray} 
\Delta_1\! & \! = \! &  \!
O \left( \left\{ \sum^{J-1}_{j=j_0 +1} \sum^{J'-1}_{j'=m'_0}+\sum^{J-1}_{j=m_0} \sum^{J'-1}_{j'=j'_0 +1} \right \} 
A^2  2^{-2js'_1 - 2j's'_2} \rkr,  \label{Delta1}\\
\Delta_2  \! & \! = \! &  \! 
O \left(  \sum^{j_0}_{j=m_0} \sum^{j'_0}_{j'=m'_0} \eps^2 \ln (1/\eps) \ 2^{j(2\nu+1)+ j'}\, 
\II \lkr 2^{j(2\nu+1)+ j'} \leq  \chia^{d-1} \rkr \rkr,  \label{Delta2}\\
\Delta_3  \! & \! = \! &  \!  
O \left(  \sum^{j_0}_{j=m_0} \sum^{j'_0}_{j'=m'_0} 
A^{p'} 2^{-p' j s'_1 - p' j' s'_2} 
\lkr \eps^2  \ln (1/\eps) 2^{2 j \nu} \rkr^{1-p'/2} \II \lkr 2^{j(2\nu+1)+ j'} > \chia^{d-1} \rkr
  \rkr, \ \ \  \label{Delta3} 
 \end{eqnarray} 
where $d$ is defined in \fr{d_value}. 
It is easy to see that for $\Delta_1$ given in \fr{Delta1} and $j_0$ and $j'_0$ given by \fr{j0j10}, as $\eps \to 0$, one has 
\be \label{Del1}
\Delta_1 = O \lkr A^2\, \chia^d \rkr,
\ee
For $\Delta_2$ defined in \fr{Delta2},  obtain 
\be \label{Del2}
\Delta_2 = O \lkr \eps^2 \ln (1/\eps) \chia^{d-1} \rkr = O \lkr A^2\,  \chia^d \rkr,  
\quad \eps \to 0.
\ee
In order to construct   upper bounds for $\Delta_3$  in \fr{Delta3}, 
 we need to consider three different cases. 
\\


\noindent
 \underline{Case 1: $s_1 \geq s_2 (2\nu +1)$. }
In this case, $d = 2 s_2/(2 s_2 +1)$ and 
\begin{eqnarray*}
\Delta_3  & \leq &  C A^2 (\chia)^{1-p'/2}\ \sum_{j={m_0}}^{j_0} 2^{-j[p' s'_1 - 2\nu( 1-p'/2)]} \ 
\sum_{j'=m'_0}^{j'_0}   2^{-  p' j' s'_2}  \, \II \lkr 2^{j'}  > (\chia)^{d-1}  2^{-j(2\nu+1)} \rkr \\
 & \leq &  C A^2 (\chia)^{(1-p'/2) + p' s'_2 (1-d)}\ 
\sum_{j={m_0}}^{j_0} 2^{-j[ p' s'_1 - 2\nu( 1-p'/2) - p'(2\nu+1)s'_2]} \\
 & = &  C A^2 (\chia)^{d} \ \sum_{j={m_0}}^{j_0} 2^{-j[p' s_1 - p' s_2(2\nu +1)]},
\end{eqnarray*}
so that, as $\eps \to 0$,  
\be \label{Del3_case1}
\Delta_3  =  O \lkr A^2\, \chia^d \, [\ln(1/\eps)]^{\II(s_1 = s_2 (2\nu +1))} \rkr.
\ee


\noindent
 \underline {Case 2: $(\frac{1}{p}-\frac{1}{2})(2\nu +1)< s_1 < s_2 (2\nu +1)$. } 
In this case, $d = 2 s_1/(2 s_1 +2\nu + 1)$ and 
\begin{eqnarray*}
\Delta_3  & \leq &  C A^2 (\chia)^{1-p'/2}\ \sum_{j={m_0}}^{j_0} 2^{-j[p' s'_1 - 2\nu( 1-p'/2)]} \ 
\sum_{j'={m'_0}}^{j'_0}   2^{-  p' j' s'_2}  \, \II \lkr 2^{j}  > (\chia)^{\frac{d-1}{2\nu +1}}  2^{- \frac{j'}{2\nu+1}} \rkr \\
 & \leq &  C A^2 (\chia)^{(1-p'/2)+  p'\frac{ (1-d)}{1+2\nu} (s_1 - (2\nu +1)(1/p'-1/2)}\ 
\sum_{j'={m'_0}}^{j'_0} 2^{- j'p'[s'_2 - s_1/(2\nu +1) + (1/2-1/p')]} \\
 & \leq &  C A^2 (\chia)^d  \ \sum_{j'={m'_0}}^{j'_0} 2^{- j'p'[s_2 - s_1/(2\nu+1)]},
\end{eqnarray*}
so that, as $\eps \to 0$,  
\be \label{Del3_case2}
\Delta_3  =  O \lkr A^2\, \chia^d  \rkr.
\ee


\noindent
 \underline {Case 3: $s_1 \leq (\frac{1}{p}-\frac{1}{2})(2\nu +1)$. } 
In this case, $d = 2 s'_1/(2 s'_1 +2\nu)$ and $p \leq 2$. Then, since 
$p  s'_1 - 2\nu( 1-p/2) = p[s_1 - (1/p - 1/2)(2\nu +1)] \leq 0$, one has
\begin{eqnarray*}
\Delta_3  & \leq &  C A^2 (\chia)^{1-p'/2}\ \sum_{j={m_0}}^{j_0} 2^{-j[p  s'_1 - 2\nu( 1-p/2)]}    \\
& \leq &  C A^2 (\chia)^{1-p'/2}\  2^{j_0 p [ (1/p - 1/2)(2\nu +1) - s_1]}\  
[\ln(1/\eps)]^{\II \lkr s_1 = (1/p - 1/2)(2\nu +1) \rkr}.
\end{eqnarray*}
Plugging in $j_0$ of the form \fr{j0j10}, obtain   as $\eps \to 0$ 
\be \label{Del3_case3}
\Delta_3  =  O \lkr A^2\, \chia^d \, [\ln(1/\eps)]^{\II( s_1 = (1/p - 1/2)(2\nu +1))} \rkr.
\ee
Now, to complete the proof, combine formulae \fr{eserto}--\fr{Del3_case3}.


\subsection{Proofs of the statements in Section \ref{sec:discrete_case}.}

{\bf Proof of Lemma \ref{lem:coef_var_dicr}. }
Subtracting $\beta_{j,k, j', k'}$  from  \fr{betaesdisver}, one obtains
\be \label{betaesdif}
\widetilde{\beta}_{j,k, j', k'} -\beta_{j,k, j', k'} =\frac{\sigma}{M}  \sum_{m \in W_j} 
\overline{\psi_{j, k, m}} \ \sum^M_{l=1}  \frac{z_m(u_l)}{g_m(u_l)}\  \eta_{j', k'}(u_l).
\ee
where $z_m(u_l)= y_m(u_l)- h_m(u_l)$. Since Fourier transform is an orthogonal transform, one has  
$\EE[ z_{m_1}(u_{l_1}) z_{m_2}(u_{l_2})]=0$ if $l_1\neq l_2$  
and $\EE[ z_{m_1}(u_{l}) z_{m_2}(u_{l})]=0$, so that   
$$
\EE[ z_{m_1}(u_{l_1}) z_{m_2}(u_{l_2})]=\frac{\sigma^2}{N} \delta(m_1-m_2)\delta(l_1-l_2).
$$
Therefore,  
\begin{eqnarray*}
\Var(\widetilde{\beta}_{j,k, j', k'}) &=&\frac{\sigma^2}{M^2N}   \sum_{m \in W_j}\left|  \psi_{j, k, m}\right|^2 
\sum^M_{l=1}\frac{1}{| g_m(u_l)|^2}|\eta_{j', k'}(u_l)|^2\\
&\asymp & 
  \frac{\sigma^2 2^{2j\nu}}{MN}\sum_{m \in W_j}\left|  \psi_{j, k, m}\right|^2 \frac{1}{M}\sum^M_{l=1}|\eta_{j', k'}(u_l)|^2 
\asymp  \frac{\sigma^2 2^{2j\nu}}{MN},
\end{eqnarray*}
 which completes the proof.
\\


{\bf Proof of Lemma \ref{lem:onedBesov}. }
Recall that
\begin{eqnarray*}
f(u, t) =  \sum_{j, k} \sum_{j', k'}   \beta_{j,k, j', k'} \psi_{j, k}(t) \eta_{j', k'}(u) \quad \mbox{and} \quad 
f_l(t)=   \sum_{j, k}    b^{(l)}_{j,k} \psi_{j, k}(t) \eta_{j', k'}(u_l),
\end{eqnarray*} 
so that 
\begin{eqnarray*}
b^{(l)}_{j,k}=   \sum^{\infty}_{j'=0} \sum_{k'\in K_l}  \beta_{j,k, j', k'}2^{j'/2}\eta(2^{j'}u_l - k'), 
\end{eqnarray*} 
where the set $K_l = \{k':\ \eta(2^{j'}u_l - k') \neq 0 \}$ is finite for any $l$ due to finite support of $\eta$.

Thus, since $p \geq 1$, for any $\delta >0$, one has
\begin{eqnarray*}
 \sumk | b^{(l)}_{j,k}|^p 
& \leq &
C  \sumk \left[ \sum^{\infty}_{j'=0} \sum_{k'\in K_l} | \beta_{j,k, j', k'}|\, 2^{j'(1+ \delta)/2}\, 2^{-j'\delta/2}\right]^2  \\
&\leq & 
C  \sumk \left( \sum^{\infty}_{j'=0} \sum_{k'\in K_l} | \beta_{j,k, j', k'}|^p\, 2^{j'(1+ \delta)p/2}\right)
\left(  \sum^{\infty}_{j'=0} \sum_{k'\in K_l}\left(2^{-j'\delta/2}\right)^{\frac{p}{p-1}}\right)^{p-1}.
\end{eqnarray*} 
Then, for any $q \geq 1$, one has
\begin{eqnarray*}
B_j= \left(  \sum^{\infty}_{j'=0}  2^{j'(1+ \delta)p/2}\,  \sum_{k, k'} | \beta_{j,k, j', k'}|^p \right)^{q/p}.
\end{eqnarray*} 
If $q/p \geq 1$, then, using Cauchy-Schwarz inequality again, it is straightforward to verify that  
\begin{eqnarray*}
B_j \leq   \tilde{C}_{\delta}  \sum^{\infty}_{j'=0} \left[ \sum_{k, k'} | \beta_{j,k, j', k'}|^p\right]^{q/p}  2^{j'(1+2 \delta)q/2}. 
\end{eqnarray*} 
Hence, 
\begin{eqnarray*}
\sum^{\infty}_{j'=0}  2^{js'_1 q} \left(  \sumk | b^{(l)}_{j,k}|^p\right)^{q/p}  
\leq \tilde{C}_{\delta}   2^{js'_1q} \ \sum^{\infty}_{j'=0} 2^{j'(1+2 \delta)q/2} \left[ \sum_{k, k'} | \beta_{j,k, j', k'}|^p\right]^{q/p} 
\leq \tilde{C}_{\delta} A^q = \tilde{A}^q
\end{eqnarray*} 
provided $s^*_2 \geq (1 + 2\delta)/2$. Since $s_2 > \max\{ 1/2, 1/p\}$ implies $s_2 >1/2$, choose $\delta = (s_2-1/2)/2$.
If $q/p <1$, then similar considerations yield
\begin{eqnarray*}
B_j \leq   \tilde{C}_{\delta}  \sum^{\infty}_{j'=0} \left[ \sum_{k, k'} | \beta_{j,k, j', k'}|^p\right]^{q/p}  2^{j'(1+ \delta)q/2}, 
\end{eqnarray*} 
so that the previous calculation holds with $\delta$ instead of $2 \delta$, and   the proof is complete.


\subsection{Proofs of the statements in Section \ref{sec:multivariate}.}

{\bf Proof of Theorem \ref{lowebdmult}. }
Repeating the proof of  Theorem \ref{lowebd} with $j'$ and $k'$ replaced by $\bjp$ and $\bkp$, respectively, and 
$s_2 j'$ replaced by $\bjp^T \bspt$, we again arrive at two cases. Denote the $r$-dimensional vector with all unit components by $\bone$.

In the  dense-dense case, we use $(r+1)$-dimensional array $w$, so that   $N = 2^{j + \bone^T \bjp}$.
Choose $\gamma^2_{j,\bjp} = A^2 2^{-j(2 s_1+1) - \bjp^T(2 \bst +\bone)}$ and observe that 
$ K(f, \tilf)  \leq   C \varepsilon^{-2}  \gamma^2_{jj'}  2^{j+\bone^T \bjp} 2^{-2\nu j}. $
Now, applying  Lemma \ref{Tsybakov} with 
\be \label{deltasqm}
\delta^2=  \gamma^2_{jj'}{2^{j+\bone^T \bjp}}/32 = A^2   2^{-  2s_1 j -  2 \bjp^T \bst}/32
\ee 
one arrives at the following   optimization problem 
\be  \label{opt1m}
(j,\bjp) = \lfi (2j s_1+ 2\bjp^T \bst) \quad \mbox{s.t} \quad   j(2s_1 + 2\nu +1) + \sum_{l=1}^r (2 s_{2,l} +1) j'_l \geq \tau_{\eps}, \  j,j'_l \geq 0\rfi,     
\ee
where $\tau_{\eps}$ is defined in formula \fr{taue}. Setting 
$j = \tau_\eps/(2 s_1 + 2\nu +1) - \sum_{l=1}^r (2 s_l +1)/(2 s_1 + 2\nu +1),$ 
arrive at optimization problem 
\be \label{optcase1m}
\bjp = \lfi \frac{2 s_1 \tau_\eps}{2 s_1 + 2\nu +1} +  \sum_{l=1}^r   \frac{2 j'_l [s_{2,l} (2 \nu +1) - s_1]}{2 s_1 + 2\nu +1}\quad  
\mbox{s.t.} \quad j'_l \geq 0, \ l=1, \cdots,r\rfi.
\ee
If $s_{2,l_0} (2 \nu +1) \geq  s_1$, then each $j'_l$ is multiplied by a nonnegative number and minimum is attained when
$j'_l =0$, $l =1, \cdots, r$. Then, $j = \tau_\eps/(2s_1 + 2\nu +1)$. On the other hand, if $s_{2,l_0} (2 \nu +1) <  s_1$,
then $j_{l_0}$ is multiplied by the smallest factor which is negative. Therefore, minimum in \fr{optcase1m} is attained
if $j=0$, $j'_l =0$, $l \neq l_0$ and $j_{l_0} = \tau_\eps/(2 s_{2, l_0} +1)$.
Plugging those values into \fr{deltasqm}, obtain
\be \label{sol:densem}
\delta^2  = \left\{ \begin{array}{ll}
C A^2\,  (\eps^2/A^2)^{\frac{2s_{2,0}}{2s_{2,0} +1}}, & \mbox{if}\ \ s_1 > s_{2,0} (2\nu +1), \\
C  A^2\,  (\eps^2/A^2)^{\frac{2s_1}{2s_1 + 2\nu + 1}}, 
& \mbox{if}\ \    s_1 \leq s_{2,0} (2\nu +1). 
\end{array} \right.
\ee

In the  sparse-dense case, we use $r$-dimensional array $w$, so that   $N = 2^{\bone^T \bjp}$.
Choose $\gamma^2_{j,\bjp} = A^2 2^{-2 j s^*_1   - \bjp^T(2 \bst +\bone)}$ and observe that 
$ K(f, \tilf)  \leq   C \varepsilon^{-2}  \gamma^2_{jj'}  2^{j+\bone^T \bjp} 2^{-2\nu j}. $
Now, applying  Lemma \ref{Tsybakov} with 
\be \label{delta2m}
\delta^2=    A^2   2^{-  2s^*_1  j -  2 \bjp^T \bst}/32
\ee 
one arrives at the following   optimization problem 
\be  \label{opt2m}
(j,\bjp) = \lfi (2j s_1+ 2\bjp^T \bst) \quad \mbox{s.t} \quad 
j(2s^*_1 + 2\nu +1) + \sum_{l=1}^r (2 s_{2,l} +1) j'_l \geq \tau_{\eps}, \  j,j'_l \geq 0\rfi,     
\ee
Again, setting 
$j = \tau_\eps/(2 s^*_1  + 2 \nu)  - \sum_{l=1}^r (2 s_l +1)/(2 s_1^* + 2 \nu),$
arrive at optimization problem 
\be \label{optcase2m}
\bjp = \lfi \frac{2 s^*_1 \tau_\eps}{2 s^*_1  + 2 \nu} +  \sum_{l=1}^r   \frac{2 j'_l [2 s_{2,l}  \nu  - s^*_1]}{2 s_1^* + 2 \nu} 
\quad  \mbox{s.t.}  \quad j'_l \geq 0, \ l=1, \cdots,r \rfi.
\ee
Repeating the reasoning applied in the dense-dense case, we obtain 
$j=0$, $j'_l =0$, $l \neq l_0$ and $j_{l_0} = \tau_\eps/(2 s_{2, l_0} +1)$ if $2 s_{2,l_0}  \nu  < s^*_1$,
and $j = \tau_\eps/(2s_1 + 2\nu +1)$, $j'_l =0$, $l =1, \cdots, r$,  if $2 s_{2,l_0}  \nu  > s^*_1$.
Plugging those values into \fr{delta2m}, obtain
\be \label{sol:sparsem}
\delta^2  = \left\{ \begin{array}{ll}
C A^2\,  (\eps^2/A^2)^{\frac{2s_{2,0}}{2s_{2,0} +1}}, & \mbox{if}\ \ 2 \nu s_{2,0} \leq s^*_1, \\
C  A^2\,  (\eps^2/A^2)^{\frac{2s^*_1}{2s^*_1 + 2\nu}}, 
& \mbox{if}\ \    2 \nu s_{2,0}  > s^*_1. 
\end{array} \right.
\ee
In order to complete the proof, combine \fr{sol:densem} and \fr{sol:sparsem} and note that $s^*_1 = s'_1$ if $p \leq 2$.
\\
 
\medskip

\noindent
{\bf Proof of Theorem \ref{uppbdm}. } Repeat the proof of Theorem \ref{uppbd} with
$j'$ and $k'$ replaced by $\bjp$ and $\bkp$, respectively,  $s_2 j'$ replaced by $\bjp^T \bspt$ and 
\bes  
2^{j_0} = (\chia)^{-\frac{d}{2 s'_1}},\quad 2^{j'_{0,l}} = (\chia)^{-\frac{d}{2 s'_{2,l}}}, \quad l=1, \cdots, r.
\ees
Then,   formulae \fr{eserto}--\fr{er2} are valid. One can also partition  $\Delta$ in \fr{delta_expr} into 
$\Delta_1$, $\Delta_2$ and $\Delta_3$ given by expressions similar to \fr{Delta1}, \fr{Delta2} and \fr{Delta3}
with $r+1$ sums in \fr{Delta1} instead of two,  $\sum^{j'_0}_{j'={m'_0}}$ replaced by $r$ respective sums and 
$\II \lkr 2^{j(2\nu+1)+ j'} > \chia^{d-1} \rkr$ replaced by $\II \lkr 2^{j(2\nu+1)+ \bone^T \bjp} > \chia^{d-1} \rkr$.
Then, upper bounds \fr{Del1} and \fr{Del2}  hold.  In order to construct   upper bounds for $\Delta_3$, 
 we again need to consider three different cases. 

In Case 1, $s_1 \geq s_{2,0} (2\nu +1)$, replace $\sum_{j'=m'_0}^{j'_0}$ by $\sum_{j'_{l_0}=m'_{l_0}}^{j'_{0,l_0}}$ 
and $\sum_{j=m_0}^{j_0}$ by the sum over $j$, $j'_1, \cdots, j'_{l_0-1}, j'_{l_0+1}, \cdots, j'_r$.
Repeating  calculations for this case, keeping  in mind that $s'_{2,l}\geq s'_{2,0}$ for any $l$
and noting that, whenever $s'_{2,l}= s'_{2,0}$, we gain an extra logarithmic factor, we arrive at 
\be \label{Del3_case1m}
\Delta_3  =  O \lkr A^2\, \chia^d \, [\ln(1/\eps)]^{\II(s_1 = s_2 (2\nu +1)) + \sum_{l\neq l_0}  \II(s_{2,l} = s_{2,0})} \rkr.
\ee
In Case 2,  $(1/p -1/2)(2\nu +1)< s_1 < s_{2,0} (2\nu +1)$, replace $\sum_{j'={m'_0}}^{j'_0}$ by 
$\sum_{\bjp \in \Upsilon (\bmp, \bjpo)}$ where $\bjpo = (j'_{0,1}, \cdots, j'_{0,r})$ and arrive at 
\fr{Del3_case2}.  In Case 3, $s_1 \leq (\frac{1}{p}-\frac{1}{2})(2\nu +1)$, since the sum over 
$\bjp$ is uniformly bounded,  calculations for the two-dimensional case hold and \fr{Del3_case3} is valid.
Combination of \fr{Del3_case1m}, \fr{Del3_case2} and \fr{Del3_case3} completes the proof.


\end{document}